\documentclass[english,11pt]{article}

\usepackage[german,french,english]{babel}
\usepackage[cp850]{inputenc}
\usepackage{latexsym,graphicx, fancybox}
\usepackage{amssymb, amsmath, amsfonts}
\usepackage{color}
\usepackage{latexsym}

\font\tenmath=msbm10 scaled 1200
\font\sevenmath=msbm7 scaled 1200
\font\fivemath=msbm5 scaled 1200

\newcommand{\vertiii}[1]{{\left\vert\kern-0.25ex\left\vert\kern-0.25ex\left\vert #1 
    \right\vert\kern-0.25ex\right\vert\kern-0.25ex\right\vert}}

\newfam\mathfam \textfont\mathfam=\tenmath
\scriptfont\mathfam=\sevenmath \scriptscriptfont\mathfam=\fivemath
\def\math{\fam\mathfam}
\def\R{{\math R}}
\def\N{{\math N}}
\def\E{{\math E}}

\def\P{{\math P}}

\newtheorem{Theorem}{Theorem}[section]

\newtheorem{Proposition}[Theorem]{Proposition}

\newtheorem{Lemma}[Theorem]{Lemma}
\newtheorem{Corollary}[Theorem]{Corollary}
\newtheorem{Remark}[Theorem]{Remark}

\newtheorem{Thm}[Theorem]{Theorem}
\newtheorem{Lem}[Theorem]{Lemma}
\newtheorem{Pro}[Theorem]{Proposition}

\newtheorem{Dfn}[Theorem]{Definition}

\newfam\mathfam \textfont\mathfam=\tenmath
\scriptfont\mathfam=\sevenmath \scriptscriptfont\mathfam=\fivemath
\def\math{\fam\mathfam}

\def \^#1{\if#1i{\accent"5E\i}\else{\accent"5E#1}\fi}

\def \ind {1 \mkern -5mu \hbox{I}}

\def \cqfd{\quad_\Box}

\begin{document}
\selectlanguage{english}
\title{\bf Optimal dual quantizers of $1D$ $\log$-concave distributions: uniqueness and Lloyd like algorithm}
 
\author{ 
{\sc Benjamin Jourdain} \thanks{Cermics, Ecole des Ponts, INRIA, Marne-la-Vall\'ee, France. E-mail: {\tt   benjamin.jourdain@enpc.fr}}
\and   
{\sc  Gilles Pag\`es} \thanks{Laboratoire de Probabilit\'es, Statistique et Mod\'elisation, UMR~8001, Campus Pierre et Marie Curie, Sorbonne Universit\'e case 158, 4, pl. Jussieu, F-75252 Paris Cedex 5, France. E-mail: {\tt  gilles.pages@upmc.fr}}}
\date{}
\maketitle 
\begin{abstract} We establish for dual quantization the counterpart of Kieffer's uniqueness result for compactly supported one dimensional probability distributions having a $\log$-concave density (also called strongly unimodal): for such distributions, $L^r$-optimal dual quantizers are unique at each level $N$, the optimal grid being the unique critical point of the quantization error.  An example of non-strongly unimodal distribution for which uniqueness of critical points fails is exhibited.
 In the quadratic $r=2$ case, we propose an algorithm to compute the unique optimal dual quantizer. It provides a counterpart of Lloyd's method~I algorithm in a Voronoi framework  (see~\cite{Lloyd,McQueen}).  Finally semi-closed forms of $L^r$-optimal dual quantizers are established for power distributions on compacts intervals and truncated exponential distributions.
  \end{abstract}


\section{Introduction} Optimal Delaunay or dual quantization has been introduced in~\cite{PaWi0} in a one dimensional setting  for probabilistic numerical   purposes, in order to produce a fast algorithm for  pricing credit derivative products in finance.  It was then developed in higher dimension in~\cite{PaWi1} as a possible  alternative to optimal Voronoi (or primal) quantization (see~\cite{IEEE, GrLu, Pag2015, PagSpring2018} for introduction)  to solve various non-linear problems in quantitative finance (American option pricing and $\delta$-hedging, stochastic control for portfolio management, etc).  Both quantization modes  are spatial discretization methods of probability distributions or random vectors, one relying on Voronoi diagrams and the other on Delaunay triangulation (in $2$-dimension). Delaunay quantization is limited to compactly supported distributions but shares a universal ``stationarity property'' (see further on) which makes it much more flexible when used as a numerical tool. This paper is essentially focused on the $1$-dimensional setting. Our  aim is to prove for Delaunay quantization some uniqueness and convergence results related to optimal  quantizers and their numerical computation  for strongly unimodal distributions known in optimal Voronoi quantization as Trushkin's and Kieffer's theorems (see e.g.~\cite{Trushkin} and~\cite{Kieff} respectively).

Let us briefly explain what Delaunay quantization is in a one dimensional setting. It answers the question: how to spatially discretize a compactly supported random variable with support (contained in) $[a,b]$ using a finite subset $\Gamma=\{x_1,\ldots,x_{N} \}\subset [a,b]$ with $x_1=a<x_2<\cdots<x_{{N-1}}<x_{N}=b$.
The basic idea to discretize a random variable $X$ is the following: when $X:(\Omega, {\cal A}, \P)\to [a,b]$ falls into the interval $[x_i,x_{i+1}]$, one replaces the value of $X$ by $\widehat X$ which take values  $x_i, x_{i+1}$  with respective probabilities  $\frac{x_{i+1}-X}{x_{i+1}-x_i}$ and $\frac{X-x_i}{x_{i+1}-x_i}$. These probabilities coming as  the coefficients when writing $X$ as of the linear interpolation of $x_i$ and $x_{i+1}$ since 
\[
X= \frac{x_{i+1}-X}{x_{i+1}-x_i}x_i + \frac{X-x_i}{x_{i+1}-x_i}x_{i+1}.
\]
 This leads to introduce the so-called {\em Delaunay} projection  (or {\em splitting operator}) defined for every $\xi\!\in [a,b]$ and $u\!\in (0,1)$ by
\begin{equation}\label{eq:ProjDel}
{\rm Proj}^{del}_{\Gamma}(\xi, u) =  \sum_{i = 1}^{N-1} \Biggl [ x_i \cdot
	\mbox{\bf 1}_{\bigl\{ 0<\, u < \frac{x_{i+1}-\xi}{x_{i+1}-  x_i} \bigr\}} + x_{i+1}\cdot \mbox{\bf 1}_{\bigl\{   \frac{x_{i+1}-\xi}{x_{i+1}-x_i} \leq u <1\bigr\}}\Biggr]
	\mbox{\bf 1}_{[x_{i}, x_{i+1})}(\xi)+b{\bf 1}_{\{\xi=b\}}
\end{equation}
  so that 
\[
\widehat X = \widehat X^{\Gamma,dual} = {\rm Proj}^{del}_{\Gamma}(X, U) \quad \mbox{ with }  \quad U\sim U\big((0,1)\big), \; U\perp \!\!\!\perp X
\]
where $\perp \!\!\! \perp$ stands for independence. The above formula can be taken as a definition for a $\Gamma$-valued dual quantizer of $X$. 

Note that, owing to the above remark that aimed at its construction 
\[
\forall\, i\!\in \{1,\ldots,N\},\quad \int_{0}^{1} {\rm Proj}^{del}_{\Gamma}(\xi, u) du = \xi 
\] 
or, equivalently,
\begin{equation}\label{eq:Statiodual}
\E\big( \widehat X\,|\, X \big) = X
\end{equation}
which is a stationarity property dual from that satisfied by quadratic optimal  primal (or {\em Voronoi}) quantization (see \eqref{eq:StatioVor} below).  In particular $X\le_{cvx}  \widehat X$ (convex ordering).
  Applications of dual quantization were first mostly devoted to provide efficient numerical schemes and fast algorithms to solve non-linear problems arising in numerical probability applied to finance like the pricing and hedging of multi-asset American style options  (see e.g.~\cite{PaWi2}) or the the pricing of credit derivatives (see~\cite{PaWi0}),  basically as a competitor of Voronoi quantization and other methods (regressions, Malliavin Monte Carlo). Its dual behaviour with respect to convex order provides an informal way to provide lower and upper-bounds in various stochastic control problems. More recently, with the development of  martingale optimal transport problems in finance, both Voronoi and Delaunay quantization methods have been shown   as a systematic tool to design time discretization schemes that preserve convex order (see~\cite{JoPa2}) and more generally to  solve   numerically discrete time martingale optimal transport problems (see~\cite{JoPa1}) which turns out to be a quite challenging problem~(see \cite{AlJo},~\cite{AlJo2},~\cite{DeMarch},~\cite{GuOb},~\cite{HL19}).


\medskip
The distribution of $\widehat X^{\Gamma,dual}$ is entirely characterized by its value set $\Gamma$ and the weights $p_i^{\Gamma} = \P(\widehat X^{\Gamma,dual}=x_i)$ given for every  $i=1,\ldots,N$, by 

\begin{align*}
 p_i^{\Gamma}&= \E\, \tfrac{X-x_{i-1}}{x_i-x_{i-1}} \mbox{\bf 1}_{\{X\in (x_{i-1},x_i]\}} +   \E\, \tfrac{x_{i+1}-X}{x_{i+1}-x_i}\mbox{\bf 1}_{\{X\in (x_i,x_{i+1}]\}}.
\end{align*}

If we introduce the cumulative distribution function  (c.d.f.) $F(x) = \P\big(X\in (-\infty,x]\big)$ and the first partial moment $K(x)= \E\, X\mbox{\bf 1}_{X\in (-\infty,x]}$, then these weights write
\begin{align}
 \nonumber  p_i^{\Gamma}&= \int_{(x_{i-1},x_i]}\tfrac{\xi-x_{i-1}}{x_i-x_{i-1}}\mu(d\xi) + \int_{(x_i,x_{i+1}]} \tfrac{x_{i+1}-\xi}{x_{i+1}-x_i} \mu(d\xi)\\
&= \frac{[K]^{x_i}_{x_{i-1}}-x_{i-1}[F]^{x_i}_{x_{i-1}} }{x_i-x_{i-1}}+\frac{[F]_{x_{i}}^{x_{i+1}} x_{i+1}-[K]_{x_{i}}^{x_{i+1}}}{x_{i+1}-x_{i}}
\label{eq:poidsdual} 
\end{align}
where, for simplicity, we will denote for a function $g:\R \to \R$ and two real numbers  $x\le y$, $[g]_x^y= g(y)-g(x)$.

The $L^r$-mean error induced by replacing $X$ by its dual quantization $\widehat X^{\Gamma, dual}$ is naturally defined by 
\begin{align*}
   \|X-\widehat X^{\Gamma, dual}\|^r_r &= \int_{a}^{b} \E\, |\xi - {\rm Proj}^{del}_x(\xi,U)|^r \mu(d\xi)\\&=\sum_{i=1}^{N-1}\int_{(x_i,x_{i+1})}\left((\xi-x_i)^r\frac{x_{i+1}-\xi}{x_{i+1}-x_{i}}+(x_{i+1}-\xi)^r\frac{\xi-x_i}{x_{i+1}-x_{i}}\right)\mu(d\xi).
\end{align*}

 \medskip
The basic application of dual quantization, like  its historical counterpart in the Voronoi sense, is to produce quadrature formulae adapted to the distribution $\mu$ of the random variable $X$ since for a  function  $g:\R\to \R$ and a $\Gamma$-quantization $\widehat X^{\Gamma}$ of $X$
\[
\E\, g(X)\simeq \E\, g\big(\widehat X^{\Gamma}\big) = \sum_{i=1}^N p_i^{\Gamma} g(x_i)
\]
where the weights $p_i^{\Gamma}$ depend on the  quantization mode (primal or dual). If  $X\!\in L^2(\P)$, and $g$ is $C^1$ and its gradient is Lipschitz continuous with constant $[\nabla g]_{\rm Lip}$ ,
writing $$g(X)-g(\widehat X)=\int_0^1\left(\nabla g(\widehat X+\alpha(X-\widehat X))-\nabla g(X)\right).(X-\widehat X)d\alpha+\nabla g(X).(X-\widehat X)$$ and using the stationarity property~\eqref{eq:Statiodual} to get rid of the expectation of the second term in the right-hand side, one obtains
\[
\big|\E\, g(X) -  \E\, g(\widehat X) \big| \le \tfrac12 [\nabla g]_{\rm Lip}  \E|X-\widehat X|^2
\]
(see~\cite{PaWi2}). Note that the counterpart of such a second order error bound in Voronoi quantization {\em only holds for optimal quadratic quantizers}. 
Nevertheless, this quadrature formula emphasizes the need for quantizers inducing, at a given level $N\ge3$,  an  as small as possible  quantization error. That is a grid $\Gamma$ such that $ \E|X-\widehat X|^2$ is minimum. This is the main purpose of optimal quantization and in fact such optimal ``minimizing'' grids do exist, see Theorem~\ref{thm:Exists} below.
Let us specify   the example of  optimal quantizations for $\mathcal{U}([0,1])$. 
It follows form~\cite{PaWi1}, Section 5.1 (see also the remark after Theorem~\ref{thm:unique}), that the $L^r$-optimal dual quantizer of  $\mathcal{U}([0,1])$ (does not depend on $r$ and) is given at level $N\ge 2$ by 
\[
\Gamma^{(N),del}  = \Big\{ \tfrac{i-1}{N-1}: i = 1, \ldots, N  \Big\}
\] 
with weights $p_1=p_N= \frac{1}{2(N-1)}$ and $p_2=p_3=\cdots=p_{N-1}=\frac{1}{N-1}$ whereas, for Voronoi quantization,  the $L^r$-optimal quantizer of  $\mathcal{U}([0,1])$ does not depend on $r$ either and is given (see~\cite{PagSpring2018}) at level $N\ge 1$ by
\[
\Gamma^{(N),vor}  = \Big\{ \tfrac{2i-1}{2N}: i = 1, \ldots, N  \Big\}
\]
with all weights given by $p_i= \frac{1}{N}$. Note that the optimal Voronoi $N$-quantizer is made up with  the midpoints of the optimal Delaunay $(N+1)$-quantizer. Consequently, in this elementary framework, Voronoi optimal $N$-quantizers correspond to midpoint quadrature formula for numerical integration over $[0,1]$ whereas  Delaunay quantization yields  the trapezoid quadrature formula. Such a property no longer holds for general distributions.

\medskip
When $X$  is an $\R^d$-valued  random vector with compactly supported distribution $\mu$, $d\ge 2$, one considers grids $\Gamma=\{x_1,\ldots,x_{N}\}\subset \R^d$ such that ${\rm supp}(\mu) \subset {\rm conv}(\Gamma)$ and the  Delaunay projection operator is defined on a  hyper-triangulation of  ${\rm conv}(\Gamma)$ sharing some minimality properties. The main feature of such dual quantization in higher dimension is that it still satisfies for every grid $\Gamma$ the above dual stationarity property.  This has been established in~\cite{PaWi1} in full generality with a natural extension to unbounded random vectors (to the price of a partial loss of the  stationarity property).
 
 Then, for any fixed $r\!\in [1, +\infty)$, one may define the lowest possible $L^r$-error induced by replacing $X$ by  any of  its dually quantization  $\widehat X^{\Gamma,dual}$ where $\Gamma$ runs over grids of size (or cardinality) at most $N$. To keep sense one should assume that $N\ge d_{\mu}+1$ where $d_{\mu}$ denotes the dimension of the vector space spanned by ${\rm supp}(\mu)$ in $\R^d$. So we define for $N\ge d_{\mu}+1$, the $L^r$ dual quantization error modulus by 
 \[
 d_{r,N}(X)  := \inf\big\{ \big\| X-\widehat X^{\Gamma, dual} \big\|_r, \; {\rm conv}(\Gamma) \supset  {\rm supp}(\mu),\; {\rm card}(\Gamma)\le N\big\}.
 \]
 It turns out (see again~\cite{PaWi1}) that it satisfies the more general bound
 \begin{align*}
 d_{r,N}(X) & =   \inf_{ Y}\Big\{ \big\|X - Y\big\|_r: Y:(\Omega\times \Omega_0, {\cal A}\otimes {\cal A}_0,\P\otimes \P_0)\to
\R^d,  \\
&  \qquad\qquad\qquad\qquad{\rm card}{(Y(\Omega\times\Omega_0))} \leq N \text{ and } \E_{\P\otimes \P_0}(Y|X) = X \Big\}.
\end{align*}
 which emphasizes the connections with martingale optimal transport explored in other papers~\cite{JoPa1, JoPa2} on the one hand and with Voronoi/primal quantization.
 
  Indeed if one replaces the above Delaunay projection by a (Borel) {\em nearest neighbour projection} on the grid $\Gamma$, denoted ${\rm Proj}_{\Gamma}^{vor}$ and if we set if $\widehat X^{\Gamma,vor}= {\rm Proj}_{\Gamma}^{vor}(X)$ for some $L^r$-integrable  random vector, then 
 \begin{align*}
e_{r,N}(X) & := \inf \big\{ \big\| X-\widehat X^{\Gamma, vor} \big\|_r,\;  {\rm card}(\Gamma)\le N\big\}  \\
& =    \inf_{Y}\Big\{ \big\|X - Y\big\|_r: Y: (\Omega , {\cal A},\P)\to\R^d,  \;{\rm card}\big(Y(\Omega)\big) \leq N   \Big\}
 \end{align*}
 One has $e_{r,n}(X) \le d_{r,N}(X)$ when both moduli make sense since dual quantization takes into account the additional martingale transport property between $X$ and its quantization. 
Note that in fact both $d_{r,N}(X)$ and $e_{r,N}(X)$ only depend on the distribution, say $\mu$, of $X$ so that we will also denote $d_{r,N}(\mu)$  (and $e_{r,N}(\mu)$).
  
  It is classical background (see~\cite{GrLu} or~\cite{Pag2015}) that the infimum is in fact a minimum and that at each {\em level } $N$ there exist an optimal grid $\Gamma^{r,vor}_{N}$ such that $e_{r,N}(X)  =  \big\| X-\widehat X^{\Gamma^{r,vor}_{N}, vor} \big\|_r$.  It should be noticed that, whereas all dual quantizations satisfy the above stationarity equation~\eqref{eq:Statiodual}, only $L^r$-optimal Voronoi quantizers with $r=2$    satisfy a stationarity property, namely a reverse one
 \begin{equation}\label{eq:StatioVor}
 \E \big(X \,|\, \widehat X^{\Gamma^{2,vor}_{N},vor}\big) = \widehat X^{\Gamma^{2,vor}_{N},vor}.
\end{equation}

Likewise, as soon as $d_{r,N}(X)<+\infty$, it is established in~\cite{PaWi1} that
  $d_{r,N}(X) $ holds as a minimum i.e. $d_{r,N}(X) =   \big\| X-\widehat X^{\Gamma^{r,dual}_{N}, dual} \big\|_r$.  To be more precise we state the original existence result for dual quantization, see~\cite{PaWi1}. 
 \begin{Thm}[Existence of optimal dual quantizers] \label{thm:Exists} Let $r\!\in [1, +\infty)$ and let $\mu $ be a compactly supported distribution on $(\R^d, {\cal B}or(\R^d))$. For every level $N\ge d_{\mu}+1$, 
 there exists at least one $L^r$-optimal grid $\Gamma^{r,del}_{_N}$ with size at most $N$ i.e. satisfying
 $$ 
 d_{r,N}(\mu)= \left(\int_{\R^d\times [0,1]} |\xi - {\rm Proj}^{del}_{\Gamma^{r,del}_{_N}}(\xi,u)|^r \mu(d\xi)du \right)^{\!\!1/r} \hskip -0,35cm = \big\|X-{\rm Proj}^{del}_{\Gamma^{r,del}_{_N}}(X,U)\big\|_r,\;  (X,U)\sim\mu\otimes U([0,1]).
 $$
Moreover ${\rm conv}\big( \Gamma^{r,del}_{_N}\big) \supset {\rm supp}(\mu)  $. If ${\rm supp}(\mu)$ has at least $N$ elements, then $\Gamma^{r,del}_{_N}$ has full size $N$ and $d_{r,N}(\mu)$ decreases to $0$ as long as it does not vanish, which never occurs  if  ${\rm supp}(\mu)$ is infinite.
 \end{Thm}
 
 

%
\providecommand{\norm}[1]{\lVert#1\rVert}
\providecommand{\pnorm}[1]{\norm{#1}_p}
\providecommand{\abs}[1]{\lvert#1\rvert}
\providecommand{\enorm}[1]{\abs{#1}}
\providecommand{\ind}[1]{\mathbbm{1}_{#1}}
\providecommand{\card}[1]{\cardop{#1}}
\newcommand{\Qn}{Q_{\text{min}}}
\newcommand{\Qx}{Q_{\text{max}}}
\newcommand{\Qkn}{Q^k_{\text{min}}}
\newcommand{\Qkx}{Q^k_{\text{max}}}
\newcommand{\qXk}{\widehat X_k}
\newcommand{\qXkk}{\widehat X_{k+1}}
\newcommand{\tp}{\pi^k_{ij}}

Finally, we recall below the main result established in~\cite{PaWi3} which is  counterpart for dual quantization of the celebrated Zador theorem ruling the sharp rate of decay to $0$ of the optimal $L^r$-quantization error and its non-asymptotic version, counterpart of Pierce's lemma. It rules the dual quantization error rate in a quite similar way for bounded random vectors.


\begin{Thm}[Rate of decay of optimal   dual quantization]\label{thm:zadoretpierce2} 
$(a)$ {\sc Sharp rate  for dual quantization}: Let $X\!\in L_{\R^d}^{\infty}(\Omega,{\cal A}, \P)$ be a bounded random vector with distribution $\P_{X}= \varphi.\lambda_d\stackrel{\perp}{+}\nu_{_X}$ where $\lambda_d$ denotes the Lebesgue measure and $\nu_{_X}$  denotes its  singular component. Then, for every $r\!\in (0, +\infty)$, 
\[
\lim_{N\to+\infty} N^{\frac 1d} d_{r,N}(X) = \widetilde J^{del}_{d,r}\left(\int_{\R^d} \varphi^{\frac{d}{d+r}}d\lambda_d\right)^{\frac 1d +\frac 1r}
\]
where $\displaystyle \widetilde J^{del}_{d,r}=\inf_{N\ge 1} N^{\frac 1d} d_{r,N}\big(\mathcal{U}([0,1]^d)\big)\ge \widetilde J^{vor}_{d,r}=\inf_{N\ge 1} N^{\frac 1d} e_{r,N}\big(\mathcal{U}([0,1]^d)\big)$.

\smallskip
When $d=1$, $\widetilde J^{del}_{1,r}=  \Big(\frac{2}{(r+1)(r+2)}\Big)^{1/r}$ whereas $\widetilde J^{vor}_{d,r}=\Big(\frac{1}{(r+1)2^r}\Big)^{1/r}$. Hence, $\frac{\widetilde J^{del}_{1,r}}{\widetilde J^{vor}_{1,r}}= \left(\frac{2^{r+1}}{r+2}\right)^{1/r} \uparrow 2$ as $r\uparrow +\infty$.

\medskip
\noindent $(b)$ {\sc Non-asymptotic bound}: Let $r, \eta >0$. For every dimension $d\ge 1$, there exists a real constant $\widetilde C^{del}_{d,\eta, r}>0$ such that, for every  random vector $X:(\Omega,{\cal A}, \P) \to \R^d$, $L^{\infty}(\P)$-bounded, 
\begin{equation}\label{eq:Zadornonasymp}
d_{r,N}(X)\le \widetilde C^{del}_{d,\eta,r} N^{-\frac 1d} \sigma_{r+\eta}(X)
\end{equation}
where, for every $p>0$, $\sigma_p(X)= \inf_{a\in \R^d}\|X-a\|_p<+\infty$.
\end{Thm}

\noindent {\bf Remark.} Note that claim~$(b)$ remains true if the support of $\P_{X}$ does not span $\R^d$ as an affine space, but $A_{\mu}$ with dimension $d'<d$. However, if  such is the case~\eqref{eq:Zadornonasymp} is suboptimal since it still  holds with  $N^{-1/d'}$ replacing $N^{-1/d}$ in the right-hand side.

%
%
 

 \medskip One of the first striking  theoretical results on optimal Voronoi quantization, beyond the existence of optimal quantizers for general distributions in any  dimension and at any level,  was obtained by Trushkin who proved (see~\cite{Trushkin}, see also~\cite{Kieff}) the following  uniqueness  result for one dimensional  strongly unimodal distributions. 
 \begin{Thm}[Trushkin, 1982] Assume that $\mu=f . \lambda$ where $f:\R\to \R_+$  is a $\log$-concave density on the real line. Let $r\ge 1$. For every integer $N\!\in \N$,  there exists a unique $L^r$-optimal grid $\Gamma_N^{r,vor}=\{x_1,\ldots,x_{N}\}\subset {\rm conv}\big({\rm supp}(\mu)\big)$ of size $N$ 
 such that for $X\sim\mu$, 
 \[
e_{r,N}(X) = \big\| X-\widehat X^{\Gamma_N^{r,vor},vor}\big\|_r.
 \]
 \end{Thm} 
 As a second step, Kieffer established in the quadratic case $r=2$,  the convergence of  the so-called Lloyd's Method~I (or Lloyd's algorithm, see~\cite{Lloyd})  at an exponential rate for strongly unimodal distributions whose $\log$-density is not piecewise affine i.e. for strongly $\log$-concave densities (see~\cite{Kieff}).

 \bigskip
 The first main result of this paper is to prove that Trushkin's uniqueness theorem remains true  for dual quantization under the same strong unimodal assumption. Then, we propose, still in $1$-dimension, a  kind of counterpart of the Lloyd's Method 1, to compute optimal quadratic dual quantizers and we prove that this algorithm converges at  an  exponential rate, uniformly in the starting point,  under a  strong unimodality property.  
 
 Finally, we also provide more specific fast algorithms to compute   $L^r$-optimal dual quantizers for two families of distributions : power distributions over a compact interval and truncated exponential distributions.

 \section{Uniqueness of optimal scalar $L^r$-dual quantizers  $r\ge1$}

 Our aim is establish uniqueness of $L^r$-optimal dual quantizers for every $r\ge 1$ under suitable assumptions on $\mu$. We exclude the trivial case when $\mu(d\xi)=\delta_x(d\xi)$ for some $x\in\R$ and the optimal grid at each level $N$ is $\{x\}$. 
 To enable dual quantization, we suppose that $\mu$ is compactly supported and denote by $a<b$ the real numbers such that $[a,b]={\rm conv}\big({\rm supp}(\mu) \big)$. For $N\ge 2$, a dual quantization grid $\Gamma$ with ${\rm card}(\Gamma)\le N$ writes $\Gamma=\{x_1,\ldots,x_N\}$ for $x_1\le x_2\le \ldots\le x_N$ satisfying $x_1\le a$ and $x_N\ge b$. Then, since $\mu(\R\setminus[a,b])=0$, for $X\sim\mu$,
 $$\|X-\widehat X^{\Gamma, dual}\|^r_r = \sum_{i=1}^{N-1}\int_{(x_i\vee a,x_{i+1}\wedge b)}\left(|\xi-x_i|^r\frac{x_{i+1}-\xi}{x_{i+1}-x_{i}}+|x_{i+1}-\xi|^r\frac{\xi-x_i}{x_{i+1}-x_{i}}\right)\mu(d\xi).$$
 When $\xi\in(x_i\vee a,x_{i+1}\wedge b)$, then, by convexity of $x\mapsto |\xi-x|^r$, one has
\begin{align*}
  &|\xi-x_i\vee a|^r\le |\xi-x_i|^r\frac{x_{i+1}-x_i\vee a}{x_{i+1}-x_i}+|x_{i+1}-\xi|^r\frac{x_i\vee a-x_i}{x_{i+1}-x_i}\\
  &|\xi-x_{i+1}\wedge b|^r\le |\xi-x_i|^r\frac{x_{i+1}-x_{i+1}\wedge b}{x_{i+1}-x_i}+|x_{i+1}-\xi|^r\frac{x_{i+1}\wedge b-x_i}{x_{i+1}-x_i}\mbox{ so that }\\
&|\xi-x_i\vee a|^r\frac{x_{i+1}\wedge b-\xi}{x_{i+1}\wedge b-x_{i}\vee a}+|x_{i+1}\wedge b-\xi|^r\frac{\xi-x_i\vee a}{x_{i+1}\wedge b-x_{i}\vee a}\le |\xi-x_i|^r\frac{x_{i+1}-\xi}{x_{i+1}-x_{i}}+|x_{i+1}-\xi|^r\frac{\xi-x_i}{x_{i+1}-x_{i}}\end{align*}
with  strict  first  inequality if $x_i<a$ (by strict convexity of $x\mapsto |\xi-x|^r$ when $r>1$ and, when $r=1$, since $\xi-x_i$ and $x_{i+1}-\xi$ have opposite signs), strict second inequality if $x_{i+1}>b$ and therefore strict third inequality if $x_i<a$ or $x_{i+1}>b$. Hence $\|X-\widehat X^{\Gamma, dual}\|^r_r$ is not smaller than
 $$\sum_{i=1}^{N-1}\int_{(x_i\vee a,x_{i+1}\wedge b)}\left(|\xi-x_i\vee a|^r\frac{x_{i+1}\wedge b-\xi}{x_{i+1}\wedge b-x_{i}\vee a}+|x_{i+1}\wedge b-\xi|^r\frac{\xi-x_i\vee a}{x_{i+1}\wedge b-x_{i}\vee a}\right)\mu(d\xi),$$
and even larger if there exists $i\in\{1,\ldots,N-1\}$ such that $x_i<a<x_{i+1}$ or $x_i<b<x_{i+1}$ (because of the support condition, $\mu$ gives positive weight to any interval $[a,x]$ and $[x,b]$ with $a<x<b$). Therefore the optimal grid for $N=2$ is $\{a,b\}$ and, when $N\ge 3$, the grid $\{a\vee x_1\wedge b,\ldots,a\vee x_N\wedge b\}$ outperforms $\Gamma$ or performs as well but contains at most $N-1$ points. If ${\rm supp}(\mu)$ contains at least $N$ points, by Theorem~\ref{thm:Exists}, any optimal grid with size at most $N$ contains $N$ points and we deduce that such a grid writes  $\{a,x_2,\ldots,x_{N-1},b\}$ with $a<x_2<\ldots< x_{N-1}<b$. Let ${\cal S}^{a,b}_{N}=\{(x_2,\ldots,x_{N-1})\in(a,b)^{N-2}:x_2<\ldots< x_{N-1}\}$ with closure $\overline{\cal S}^{a,b}_{N}=\{(x_2,\ldots,x_{N-1})\in[a,b]^{N-2}:x_2\le \ldots\le x_{N-1}\}$ and for $x=(x_2,\ldots,x_{N-1})\in{\cal S}^{a,b}_{N}$, \[
L_{_{\!N}}(x):=\|X-\widehat X^{\Gamma, dual}\|^r_r \mbox{ where }\Gamma=\{a,x_2,\ldots,x_{N-1},b\}.
\] 
We will of course use the natural convention $x_1=a$ and $x_N=b$ in what follows. The optimal grids, which exist according to Theorem~\ref{thm:Exists}, are of the form $\{a,x_2,\ldots,x_{N-1},b\}$ with $x=(x_2,\ldots,x_{N-1})\in{\cal S}^{a,b}_{N}$ minimizing $L_N$ over ${\cal S}^{a,b}_{N}$ when ${\rm supp}(\mu)$ contains at least $N$ points.
\begin{Thm}[Uniqueness of critical points of $L_N$] \label{thm:unique} Let $N\ge 3$ and $r\!\in [1, +\infty)$. Assume that $\mu([a,b])=1$ and $\mu$ is atomless. Then the function $L_{_{\!N}} : {\cal S}^{a,b}_{N}\to \R_+$ defined just above is differentiable and its gradient $\nabla L_N$ admits a continuous extension on $\overline{{\cal S}}^{a,b}_{N}$. If $\mu$ admits a density $f$ with respect to the Lebesgue measure which is positive and $\log$-concave on the interval $(a,b)$ and vanishes outside
  , then the continuous extension of $\nabla L_N$ has a unique  zero $x^{\star}$ and $x^{\star}\in{\cal S}^{a,b}_{N}$.
\end{Thm}
We deduce the following result.
\begin{Corollary}[Uniqueness of $L^r$-quantizers]\label{cor:unique} 
  Let $N\ge 3$ and $r\!\in [1, +\infty)$. Assume that $\mu$ is atomless and ${\rm conv}\big({\rm supp}(\mu) \big)=[a,b]$ with $-\infty<a<b<+\infty$. Then the $L^r$-optimal dual quantization grids at the level $N$ are of the form $\{a,x_2,\ldots,x_{N-1},b\}$ with $(x_2,\ldots,x_{N-1})\in{\cal S}^{a,b}_{N}$ solving the master equation $\nabla L_N(x_2,\ldots,x_{N-1})=0$. If, moreover, $\mu$ admits a density $f$ with respect to the Lebesgue measure which is positive and $\log$-concave on the interval $(a,b)$ and vanishes outside, then the unique $L^r$-optimal dual quantization grid at level $N$ of $\mu$ is $\{a,x^{\star}_2,\ldots,x^{\star}_{N-1},b\}$ where $x^{\star}$ is the unique critical point of $L_N$ in ${\cal S}^{a,b}_{N}$.
\end{Corollary}
Notice that by concavity of $x\mapsto \log f(x)$ on $(a,b)$, this function is continuous on $(a,b)$ and admits limits in $\{-\infty\}\cup\R$ as $x\to a+$ and $x\to b-$ so that $f$ is continuous on $(a,b)$ and admits limits in $\R_+$ as $x\to a+$ and $x\to b-$. The density $f$ is continuous on $\R$ iff both limits are equal to $0$. In any case, $f$ is bounded on $\R$. 

To prove Theorem~\ref{thm:unique}, we will rely on the following two classical results, tailored  variants of the celebrated Mountain pass lemma and Gershgorin's lemma respectively.

\begin{Thm}[Mountain pass Lemma]\label{thm:MPL}  

{\bf  Compact case} (see \cite{LaPa}) Assume that $K\!\subset \R^n$ is the closure of a
nonempty compact convex open  set O (then $\stackrel{\circ}{K} = O$) and that $L:K\longrightarrow \R$ is $C^1$
on
$\stackrel{\circ}{K}$, $\nabla L$ admits a continuous extension  on $K$
satisfying $\{\nabla L\!=\!0\}\!\subset \,\stackrel{\circ}{K}$ and for every small enough $
\varepsilon>0$, $(Id-\varepsilon \nabla L)(\stackrel{\circ}{K})\!\subset \stackrel{\circ}{K}$. If two distinct zeros of $\nabla L$ are (strict) local
 minima, then $\nabla L$ has  a third zero  which can in no case be a local minimum.
\end{Thm}

\begin{Proposition}[\`A la Gershgorin Lemma] \label{lem:Gershg+}$(a)$ Let $A= [a_{ij}]$ be an $n\times n$ symmetric matrix with  dominating diagonal i.e.
\[
\forall \, i\!\in \{1,\ldots,n\},\;  \forall\, j \!\in\{1, \ldots,n  \}\setminus\{i\} , \; a_{ij} \le 0\; \mbox{ and }\;  \Lambda_i := \sum_{j=1}^n a_{ij}\ge 0.
\]
Then all eigenvalues of $A$ are non-negative. 

\smallskip 
\noindent $(b)$ If moreover,  $A$ is tridiagonal with $a_{ii\pm1}<0$, for  $i=2,\ldots,n-1$, $a_{12}$, $a_{nn-1}<0$, $\Lambda_1$ or $\Lambda_{n}>0$, then all eigenvalues of $A$ are positive.
\end{Proposition}

\noindent {\bf Proof.} Let $\lambda\!\in\R$ be an eigenvalue of the symmetric matrix $A$ and $x=(x_1,\ldots,x_n)$ one of its eigenvectors. Let $i_0\in{\rm argmax}_{1\le i\le n} |x_i|$. We may assume without loss of generality that
$x_{i_0}=1$. Also note that $a_{ii}\ge 0$, $i=1,\ldots,n$.

$(a)$    We have $|x_j|\le 1$ for all $j=1,\ldots,n$ so that 
\[
\lambda =a_{i_0i_0} + \sum_{j\neq i_0}a_{i_0j} x_j\ge a_{i_0i_0} + \sum_{j\neq i_0}a_{i_0j} |x_j|\ge \Lambda_{i_0}\ge 0.
\]

\noindent $(b)$  If $\lambda =0$, then under the convention $a_{10}=a_{nn+1}=0$ and $x_0=x_{n+1}=1$, we have
$$  
0= \lambda =a_{i_0i_0}+ a_{i_0i_0-1} x_{i_0-1}+ a_{i_0i_0+1} x_{i_0+1} \ge  a_{i_0i_0}+ a_{i_0i_0-1} |x_{i_0-1}| + a_{i_0i_0+1} |x_{i_0+1}| \ge  \Lambda_i =0
$$ 
so that $x_{i_0\pm1}= |x_{i_0\pm1}| =1$. Then, by induction, we show that $x_{i}=1$ for all $i=1,\ldots,n$ i.e. $x= \mbox{\bf 1}$. But then, if $L_1=a_{11}+a_{12}>0$, 
$$
\lambda =\lambda x_1 =  a_{11}x_1+a_{12}x_2 = a_{11}+a_{12}>0
$$ 
which yields a contradiction. One concludes likewise if $\Lambda_{n} >0$. Hence $\lambda>0$.\hfill$\Box$

\bigskip
In order to prove Theorem~\ref{thm:unique} we will follow the strategy originally developed in~\cite{LaPa} for Voronoi optimal quantizers.

\medskip
\noindent {\bf Proof.} In the proof, the hypotheses on $\mu$ will only be gradually reinforced to those made in the statement. 

\begin{figure}[h!]
 \centering  
  \begin{minipage}[h]{0.2\textwidth}
   \centering
   \includegraphics[height= 2cm, width=3cm]{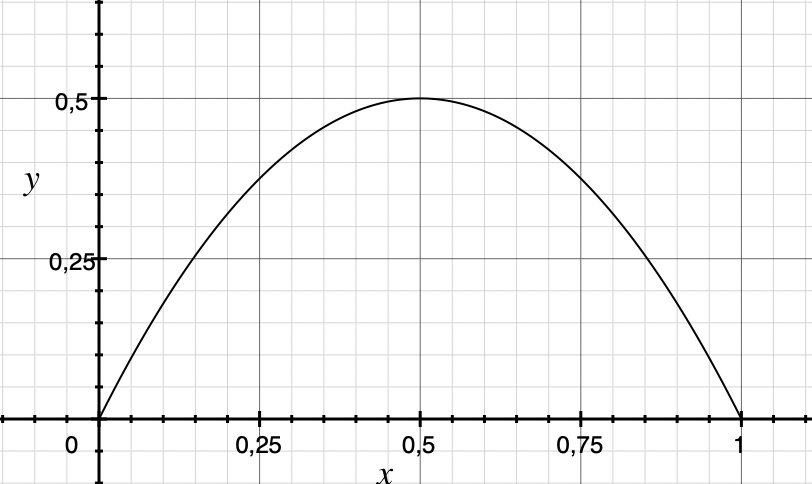} 
   \end{minipage}
 \begin{minipage}[h]{0.2\textwidth}
 \centering
 \includegraphics[height= 2cm, width=3cm]{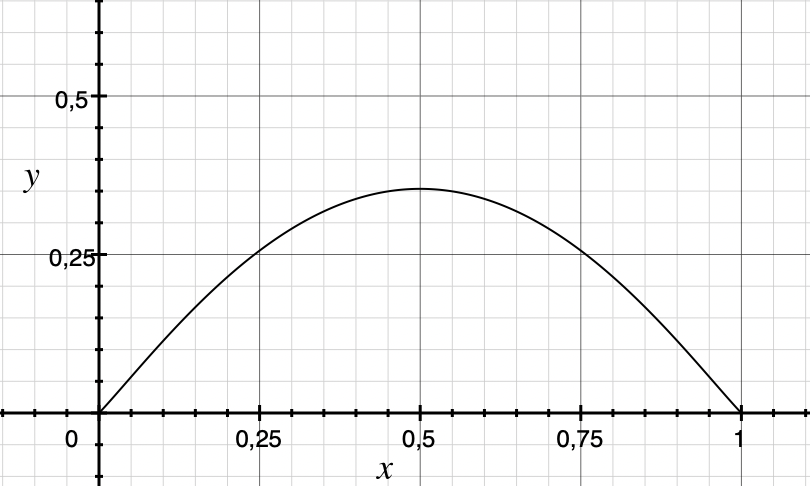}
 \end{minipage}
 \begin{minipage}[h]{0.2\textwidth}
   \centering
    \includegraphics[height= 2cm, width=3cm]{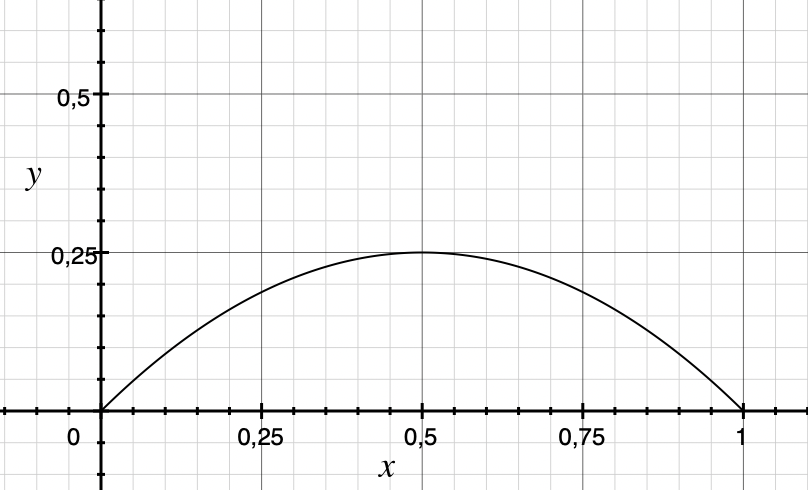}
   \end{minipage}
   \centering  
  \begin{minipage}[h]{0.2\textwidth}
   \centering
   \includegraphics[height= 1cm, width=3cm]{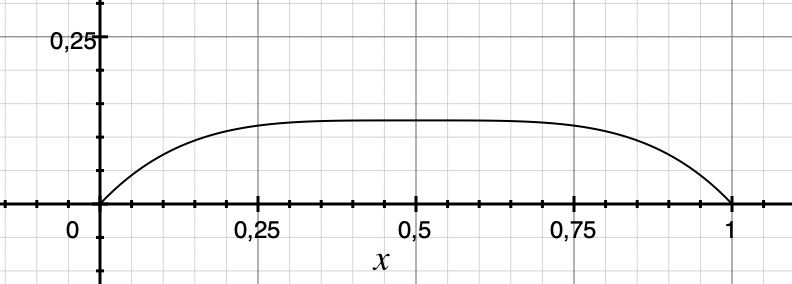} 
   \end{minipage}
 \vskip 0,5cm 
 \begin{minipage}[h]{0.2\textwidth}
 \centering
 \includegraphics[height= 1cm, width=3cm]{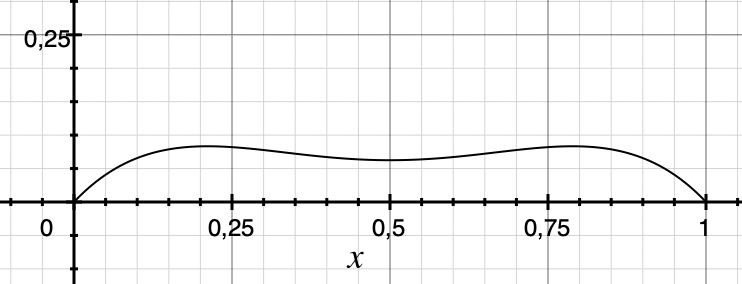}
 \end{minipage}
 \begin{minipage}[h]{0.2\textwidth}
   \centering
    \includegraphics[height= 1cm, width=3cm]{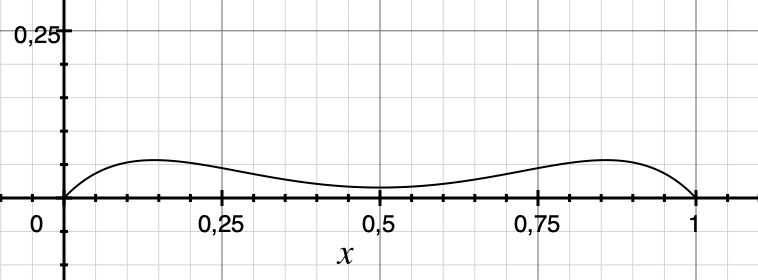}
   \end{minipage}
 \begin{minipage}[h]{0.2\textwidth}
   \centering
    \includegraphics[height= 1cm, width=3cm]{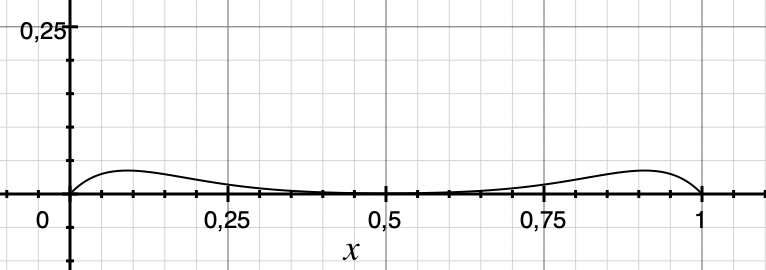}
   \end{minipage}

 \caption{\it Functions $\varpi_r$ for (from left to right) $r=1,1.5,2,3,4,6,10.$}
     \label{fig:varphi_r}
  \end{figure}

We first assume that $\mu([a,b])=1$. We  know that, for every $\xi \!\in [x_{i-1},x_i]$ and every $u \!\in   [0,1]$,
\[
{\rm Proj}^{del}_x(\xi,u) =\sum_{i=2}^{N}  \mbox{\bf 1}_{\{0\le u < \frac{\xi-x_{i-1}}{\Delta x_i}  \}}x_i  + \mbox{\bf 1}_{\{  \frac{\xi-x_{i-1}}{\Delta x_i} \le u \le 1 \}}x_{i-1}, \qquad \Delta x_i = x_i-x_{i-1}, \; i=2,\ldots,N.
\]
Now 
\begin{align}
\nonumber L_{_{\!N}}(x) &=\int_{x_{i-1}}^{x_i} \Big[ |\xi-x_{i-1} |^r   \frac{x_{i}-\xi}{\Delta x_i}+|x_i-\xi|^r \frac{\xi - x_{i-1}}{\Delta x_i}   \Big] \mu(d\xi) \\
\label{eq:LrDistor1D}
 &= \sum_{i=2}^{N} (\Delta x_i)^{r} \int_{x_{i-1}}^{x_i}\varpi_r\Big( \frac{\xi-x_{i-1}}{\Delta x_i} \Big) \mu(d\xi)
\end{align}
where 
$$
\varpi_r(u) = u^{r}(1-u) + (1-u)^r u,\quad u \!\in [0,1],
$$
(See Figure~\ref{fig:varphi_r}). Note that $\varpi_r\Big( \frac{x_i-x_{i-1}}{\Delta x_i} \Big)\mu(\{x_i\})=0=\varpi_r(\frac{x_{i-1}-x_{i-1}}{\Delta x_i})\mu(\{x_{i-1}\})$ so that the notation $\int_{x_{i-1}}^{x_i}\varpi_r\Big( \frac{\xi-x_{i-1}}{\Delta x_i} \Big) \mu(d\xi)$ makes sense even if $\mu$ weights points. The function $\varpi_r$ when extended under the same notation  by the value $0$ outside the interval $[0,1]$ is continuous and bounded by $1/2$ on the real line. As a consequence, Lebesgue's theorem ensures that $(y,z)\mapsto \int_y^z \varpi_r(\frac{\xi-y}{z-y})\mu(d\xi)=\int_\R\varpi_r(\frac{\xi-y}{z-y})\mu(d\xi)$ is continuous on $\{(y,z)\in\R^2: y<z\}$. Moreover, for $y<z$, $0\le (z-y)^r\int_y^z \varpi_r(\frac{\xi-y}{z-y})\mu(d\xi)\le \frac{(z-y)^r}{2}$. We deduce that $L_N$ is continuous on $ {\cal S}^{a,b}_{N}$ and can be continuously extended to its closure $K= \overline{{\cal S}}^{a,b}_{N}$.

When $\mu$ has a density $f$, an elementary change of variable in each integral yields the alternative formulation
\begin{equation}\label{eq:LrDistor1Db}
L_{_{\!N}}(x)= \sum_{i=2}^{N} (\Delta x_i)^{r-1}  \int_{0}^{1} \varpi_r(z)   f\big(x_{i-1} +z\Delta x_i\big) dz.
\end{equation}

Since the extended function $\varpi_r$ is differentiable outside $\{0,1\}$ with a bounded derivative, Lebesgue's theorem ensures that $(y,z)\mapsto \int_\R \varpi_r(\frac{\xi-y}{z-y})\mu(d\xi)$ admits a partial derivative with respect to its first (resp. second) variable equal to $\int_\R \frac{\xi-z}{(z-y)^2}\varpi'_r(\frac{\xi-y}{z-y})\mu(d\xi)$ (resp. $\int_\R \frac{y-\xi}{(z-y)^2}\varpi'_r(\frac{\xi-y}{z-y})\mu(d\xi)$) at each point $(y,z)$ such that $y<z$ and $\mu(\{y\})=0$ (resp. $\mu(\{z\})=0$). Hence, for $ i=2:N-1$ and $x\in{\cal S}^{a,b}_{N}$ such that $\mu(\{x_i\})=0$,  $L_N$ admits a partial derivative with respect to $x_i$ at $x$ given by
\begin{align}\label{eq:grad1}
  \partial_{x_i} L_{_{\!N}}(x) &=  (\Delta x_i)^{r-1} \int_{(x_{i-1},x_i]}\Psi_r\Big( \frac{\xi-x_{i-1}}{\Delta x_i} \Big) \mu(d\xi)- (\Delta x_{i+1})^{r-1}  \int_{[x_{i},{x_{i+1}})}\Psi_r\Big( \frac{x_{i+1}-\xi}{\Delta x_{i+1}} \Big) \mu(d\xi),\end{align}
where the function $\Psi_r$ is defined by 
\begin{equation}\label{eq:Psi_r}
\forall u\in(0,1),\;\Psi_r(u) = r \varpi_r(u)-u\varpi'_r(u) = (r-1)u(1-u)^r +u^{r+1} +ru^2(1-u)^{r-1}>0
\end{equation}
and by $\Psi_r(0)=0$ and $\Psi_r(1)=1+{\mathbf 1}_{\{r=1\}}$. 

By Fubini's theorem and since $\Psi_r(0)=0$,
\begin{align*}
   \int_{(x_{i-1},x_i]}\Psi_r\Big( \frac{\xi-x_{i-1}}{\Delta x_i} \Big) \mu(d\xi)&=\int_{(x_{i-1},x_i]}\int_0^1{\mathbf 1}_{\{z<\frac{\xi-x_{i-1}}{\Delta x_i}\}}\Psi'_r(z)dz\mu(d\xi)\\&=\int_0^1\Psi'_r(z)\big(F(x_i)-F(x_{i-1} +z\Delta x_i)\big) dz.
\end{align*}
Hence, dealing in a similar way with the second term in the right-hand side of~\eqref{eq:grad1}, we obtain the following second form of $\partial_{x_i} L_{_{\!N}}(x)$ :
\begin{align}\label{eq:grad3}
\nonumber \partial_{x_i} L_{_{\!N}}(x) &= (\Delta x_i)^{r-1} \int_{0}^{1}\Psi'_r(z)\big(F(x_i)-F(x_{i-1} +z\Delta x_i)\big) dz \\
& \quad -(\Delta x_{i+1})^{r-1}  \int_0^{1}\Psi'_r(1-z) \big( F(x_{i} +z\Delta x_{i+1})-F(x_i)\big) dz.
\end{align}
When $\mu$ is atomless, then $L_N$ is differentiable on ${\cal S}^{a,b}_{N}$ and we easily deduce from the continuity of $F$ and \eqref{eq:grad3} that $\nabla L_N$ is continuous on ${\cal S}^{a,b}_{N}$ and admits a continuous extension on $\overline{\cal S}^{a,b}_{N}$.

When $\mu$ has a density $f$, then the partial derivative is also equal to
\begin{align}\label{eq:grad2}
 \partial_{x_i} L_{_{\!N}}(x) &= (\Delta x_i)^{r} \int_{0}^{1}\Psi_r(z) f\big(x_{i-1} +z\Delta x_i\big) dz -(\Delta x_{i+1})^{r}  \int_0^{1}\Psi_r(1-z) f(x_{i} +z\Delta x_{i+1}\big) dz.
\end{align}


Note that each of these  three forms of partial derivatives $\partial_{x_i}L_{_{\!N}}(x)$ yields a version of the {\em master  equation} for dual quantization at level $N$,  $\nabla L_{_{\!N}}(x)=0$.

From now on, we assume that $\mu$ has a density $f$ which is $du$ a.e. positive on $(a,b)$. If $x\in\overline{{\cal S}}^{a,b}_{N}$ solves the master equation then by \eqref{eq:grad2}, $(\Delta x_i)^r=0\Leftrightarrow (\Delta x_{i+1})^r=0$ for $i=2,\ldots,N-1$ and since $x_1=a<b=x_N$, necessarily $\Delta x_i>0$ for $i=2,\ldots,N$ i.e. $x\in {\cal S}^{a,b}_{N}= \stackrel{\circ}{K}$.


\smallskip When the density $f$ is continuous on $(a,b)$, the cumulative distribution function $F$ is continuously differentiable on this interval and it follows from~\eqref{eq:grad3} that the Hessian of $L_N$ does exist on ${\cal S}^{a,b}_{N}$ and has a symmetric tridiagonal structure. However in order to apply the refined Gershgorin Lemma (Lemma~\ref{lem:Gershg+}$(b)$ with $n= N-2$), we need to show that the sub-diagonal terms are non-positive and  the sum of its lines i.e. $ \sum_{\ell= 0, \pm1}  \partial^2_{{x_i}x_{i+\ell}} L_{_{\!N}}(x),\;  i=2:N-1$ (with the obvious convention that $\partial_{x_1}[\ldots] =\partial_{x_{{N}}}[ \ldots]=0$). We assume that the density $f$ is also continuous and rely on~\eqref{eq:grad3} to derive that  
\begin{align}
\nonumber\partial^2_{x_{i-1}x_i} L_{_{\!N}}(x) &=  -(r-1)(\Delta x_i)^{r-2} \int_{0}^{1}\Psi'_r(z)\big(F(x_i)-F(x_{i-1} +z\Delta x_i)\big) dz\\
\nonumber& \quad - (\Delta x_i)^{r-1}\int_{0}^{1}\Psi'_r(z)(1-z) f(x_{i-1} +z\Delta x_i)dz\\
\label{eq:Hess_i-1i} & = - (\Delta x_i)^{r-1}\int_0^1 \vartheta_r(z) f(x_{i-1} +z\Delta x_i)  dz
\end{align}
where, for every $z\!\in (0,1)$,  
\[
\vartheta_r(z) = (r-1)\Psi_r(z)+(1-z) \Psi'_r(z) = (r+1)\big(z(1-z)^r +z^r(1-z)  \big) +(r-1) \big(z^{r+1}+(1-z)^{r+1}\big) >0.
\]
Hence, $f$ being $du$-a.e. positive on $(a,b)$, $\partial^2_{x_{i-1}x_i} L_{_{\!N}}(x) <0$.
 
\smallskip
Similar computations show, still assuming that $f$ is continuous,  that 
\begin{align}
\nonumber \partial^2_{x^2_i} L_{_{\!N}}(x)&= \Psi_r(1) f(x_i) \big((\Delta x_i)^{r-1}+(\Delta x_{i+1})^{r-1} \big) \\
\label{eq:Hess_ii}  & \quad+ (\Delta x_i)^{r-1}\int_0^1 \widetilde \vartheta_r(z)  f\big(x_{i-1} +z\Delta x_i\big) dz+ (\Delta x_{i+1})^{r-1}  \int_0^1 \widetilde \vartheta_r(1-z)  f\big(x_{i} +z\Delta x_{i+1}\big) dz
\end{align}
($\Psi_r(1)=1+{\bf 1}_{\{r=1\}}$) where
\[
\widetilde \vartheta_r(z)= (r-1)\Psi_r(z)-z \Psi'_r(z).
\]
Let us introduce for every $i=2,\ldots,N-1$, the quantity
\begin{align*}
      S_i &=   \Psi_r(1) \big((\Delta x_i)^{r-1} +  (\Delta x_{i+1})^{r-1} \big) f(x_i) \\
       &\hskip 2cm - (\Delta x_i)^{r-1}\int_{0}^{1}\Psi'_r(z) f(x_{i-1}+z\Delta x_i) dz- (\Delta x_{i+1})^{r-1}   \int_{0}^{1}\Psi'_r(1-z)f(x_{i}+z\Delta x_{i+1})dz.
\end{align*}
One derives from~\eqref{eq:Hess_i-1i},~\eqref{eq:Hess_ii}  and the obvious fact   $\widetilde \vartheta _r -\vartheta_r = -\Psi'_r$ that
 \[
 S_i  = \sum_{\ell= 0, \pm1}  \partial^2_{{x_i}x_{i+\ell}} L_{_{\!N}}(x)=\sum_{j=2}^{N-1}  \partial^2_{{x_i}x_{j}} L_{_{\!N}}(x) \quad \mbox{ for }\quad  i=3:N-2.
\]
(We could have taken  advantage of the anti-symmetries induced by the fact that $\partial_{x_i}\Delta x_i + \partial_{x_{i-1}} \Delta x_i =0$ to compute $S_i$ without computing $\partial^2_{x^2_i} L_{_{\!N}}(x) $ but we will need a closed form of the diagonal term of the Hessian for the counterexample below).

Moreover, one checks  that by positivity of $\vartheta_r$ on $(0,1)$ and of $f$ on $(a,b)$, 
\begin{align*}
   &\sum_{j=2}^{N-1}  \partial^2_{{x_2}x_{j}} L_{_{\!N}}(x)=\sum_{\ell=  0,1}  \partial^2_{{x_2}x_{i+\ell}} L_{_{\!N}}(x)=S_2+(\Delta x_2)^{r-1}\int_0^1 \vartheta_r(z) f(x_{1} +z\Delta x_2)  dz>S_2\mbox{ and }\\&\sum_{j=2}^{N-1}  \partial^2_{{x_{N-1}}x_{j}} L_{_{\!N}}(x)=\sum_{\ell= -1, 0}  \partial^2_{{x_{{N-1}}}x_{{N-1+\ell}}} L_{_{\!N}}(x)=S_{{N-1}}+(\Delta x_{N})^{r-1}\int_0^1 \vartheta_r(z) f(x_{N-1} +z\Delta x_{N})  dz > S_{{N-1}}.
\end{align*}
Assume now that $f$ is  positive and right differentiable on $(a,b)$ with right derivative $f'_r$. Then, by an integration by part, one shows that
\begin{align*}
 S_i &=   (\Delta x_i)^{r}\int_{0}^{1}\Psi_r(z)f'_r(x_{i-1}+z\Delta x_i) dz - (\Delta x_{i+1})^{r}   \int_{0}^{1}\Psi_r(1-z)f'_r(x_{i}+z\Delta x_{i+1}) dz \\
 &= (\Delta x_i)^{r} \hskip-0,15cm  \int_{0}^{1} \hskip-0,25cm  \Psi_r(z)  \frac{f'_r(x_{i-1}+z\Delta x_i)}{f(x_{i-1}+z\Delta x_i)} f(x_{i-1}+z\Delta x_i)dz \\
 &\qquad - (\Delta x_{i+1})^{r} \hskip-0,15cm \int_{0}^{1} \hskip-0,25cm  \Psi_r(1-z) \frac{f'_r(x_{i}+z\Delta x_{i+1})}{f(x_{i}+z\Delta x_{i+1})} f(x_{i}+z\Delta x_{i+1})dz.
\end{align*}

Now if, furthermore, $f$ is $\log$-concave  then $f$  is right differentiable and $\frac{f'_r}{f}$ is non-increasing so that $\frac{f'_r}{f}\le \frac{f'_r}{f}(x_i)$ on $(x_i,x_{i+1})$ and $\frac{f'_r}{f}\ge \frac{f'_r}{f}(x_i)$ on $(x_{i-1},x_i)$. Since $\Psi_r$ is positive on $(0,1)$, we deduce that
\begin{align*}
  S_i &\ge \frac{f'_r}{f}(x_i)\left((\Delta x_i)^{r} \hskip-0,15cm  \int_{0}^{1} \hskip-0,25cm  \Psi_r(z)  f(x_{i-1}+z\Delta x_i)dz  - (\Delta x_{i+1})^{r} \hskip-0,15cm \int_{0}^{1} \hskip-0,25cm  \Psi_r(1-z)  f(x_{i}+z\Delta x_{i+1})dz\right)
\end{align*}
where the second factor in the right-hand side is equal to $0$ when $x\in{\cal S}^{a,b}_{N}$ is solution to the master equation (derived from)~\eqref{eq:grad2}. Consequently, $\sum_{j=2}^{N-1}\partial^2_{x_jx_i}L_N(x)= S_i \ge 0$ for $i=3:N-2$ and $\sum_{j=2}^{N-1}\partial^2_{x_jx_i}L_N(x)> S_i\ge 0$ for $i\in\{2,N-1\}$. It follows  from Proposition~\ref{lem:Gershg+} that the Hessian $\nabla^2L_N(x)$ at  an equilibrium point $x\in{\cal S}^{a,b}_{N}$ has a strictly positive spectrum and $x$ is consequently a strict local minimum of $L_{_{\!N}}$ on ${\cal S}^{a,b}_{N}$. If we can prove that $(I_d-\varepsilon \nabla L)({{\cal S}^{a,b}_{N}})\!\subset {{\cal S}^{a,b}_{N}}$ for small enough $\varepsilon$, then we may apply  (the  variant of) the Mountain Pass Lemma (Theorem~\ref{thm:MPL})  to the convex compact $\overline{\cal S}^{a,b}_{N}$ of $\R^{N-2}$ with non empty interior to conclude that $L_N$ admits at most one equilibrium point $x\in{\cal S}^{a,b}_{N}$. This is the purpose of the next lemma, the hypothesis of which is satisfied when the density $f$ is positive and $\log$-concave on $(a,b)$ since, according to the remark just after Corollary~\ref{cor:unique}, $f$ is then bounded.

\begin{Lem} Let $r\!\in [1, +\infty)$. If the density $f$ satisfies 
\[
f \mbox{ bounded if $r=1$ or } \int_a^b f^{\frac{1}{r-1}}(\xi)d\xi <+\infty\mbox{ if } \; r\!\in (1,2),
\]
then, for $\varepsilon >0$ small enough,  $(Id-\varepsilon \nabla L_N)({\cal S}^{a,b}_{N})\!\subset {\cal S}^{a,b}_{N}$.
\end{Lem}

\noindent {\bf Proof.}
Let $x\in{\cal S}^{a,b}_{N}$.

\smallskip Assume $r\!\in (1,2)$. Then, for every $u,v\!\in [a,b]$, H\"older's inequality implies 
\[
\int_u^v f(\xi)d\xi \le \left(\int_a^b f^{\frac{1}{r-1}}(\xi)d\xi\right)^{r-1}   \, (v-u)^{2-r} .
\]
On the other hand, for $i=2,\ldots,N-2$,  it follows from~\eqref{eq:grad2} that 
\begin{align*}
  \partial_{x_{i+1}}L_{_{\!N}}(x)  - \partial_{x_i}L_{_{\!N}}(x) & \le   (\Delta x_{i+1})^{r}  \int_0^{1}\big(\Psi_r(z) +\Psi_r(1-z)\big)f(x_{i} +z\Delta x_{i+1}\big) dz\\
  &\le \Delta x_{i+1} \,C_r  (\Delta x_{i+1})^{r-2}  \int_{x_{i}}^{x_{i+1}} f(\xi) d\xi\\
&\le   \Delta x_{i+1}  \,C_r \left(\int_a^b f^{\frac{1}{r-1}}(\xi)d\xi\right)^{r-1}
\end{align*} 
where $C_r = \sup_{z \in [0,1]}\big(\Psi_r(z)+\Psi_r(1-z)\big)< + \infty$ according to~\eqref{eq:Psi_r}. Consequently for $\varepsilon <\Big(  \,C_r \left(\int_a^b f^{\frac{1}{r-1}}(\xi)d\xi\right)^{r-1}\Big)^{-1}$
\[
  x_i- \varepsilon  \partial_{x_i}L_{_{\!N}}(x) < x_{i+1}-\varepsilon  \partial_{x_{i+1}}L_{_{\!N}}(x), \quad i=2:N-1.
\]

\smallskip If $r=1$, this inequality follows likewise by replacing $\left(\int_a^b f^{\frac{1}{r-1}}(\xi)d\xi\right)^{r-1} $ by $\|f\|_{\infty}$.

\smallskip If $r\ge 2$, just note that $ (\Delta x_{i+1})^{r-2}  \int_{x_{i}}^{x_{i+1}} f(\xi) d\xi\le (b-a)^{r-2}$ and choose $\varepsilon \le \big(C_r(b-a)^{r-2}  \big)^{-1}$. 

 It remains to prove that $x_2-\varepsilon \,\partial _{x_2}L(x)>a$ and $x_{{N-1}}-\varepsilon\, \partial _{x_{{N-1}}}L(x)<b$. In fact 
\begin{align*}
x_2-\varepsilon \,\partial_{x_2} L(x) &> x_2-\varepsilon (\Delta x_2)^{r}\int_0^1\Psi_r(z) f(a+ z\Delta x_2)dz\\
&\ge x_2-\Delta x_2 \varepsilon \|\Psi_r\|_{\sup} (\Delta x_2)^{r-2}\int_{a}^{x_2}f(\xi)d\xi.
\end{align*} 
Inspecting the same cases as above, one shows  under the assumptions made on $f$ for $r\!\in [1,2]$, that for $\varepsilon \!\in (0, \varepsilon'_r]$ small enough (independently of $x$), 
\[
x_2-\varepsilon\, \partial_{x_2} L(x) > x_2- (x_2-a) = a.
\]
The second inequality follows likewise. This completes the proof of the lemma. \hfill $\Box$

\bigskip
The following counterexample shows that uniqueness of critical points of $L_{_{\!N}}$ may fail when the density $f$ is continuous, $du$-$a.e.$ positive and (left and) right differentiable, but not not $\log$-concave.

\bigskip

\noindent {\bf  Counter--example.}  The idea to devise this counter-example is to find a distribution $\mu$ with a periodic density on the interval $[0,1]$  that  trivially makes $x^{\star}=\big(\frac{k-1}{N-1} \big)_{k=1,\ldots,N}$  an equilibrium at level $N$ but which assigns much mass in between the codewords $x^*_k= \frac{k-1}{N-1}$  so that  this equilibrium cannot be a local minimum of $L_{_{\!N}}$. As a consequence there will be at least one further equilibrium point: the minimum of $L_{_{\!N}}$ known to lie in $S_{_N}^{0,1}$.

\begin{figure}[h!]
\centering  
   \includegraphics[height= 3cm, width=5cm]{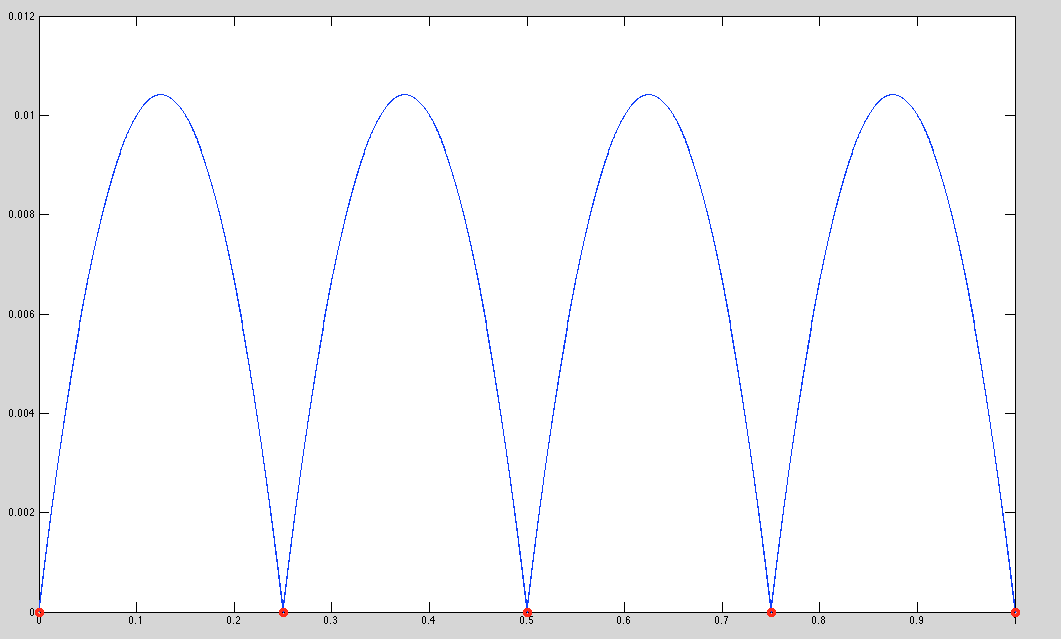}
   \caption{\em The  probability density $g$ for $r=2$ and $N=5$ (see~\eqref{eq:Ctrxple}). Red bulllets are the $5$ codewords $x^*_k$, $k=1:5$.} 
\end{figure}
 \smallskip
 Let $r\ge 1$ and let $g\!\in {\cal C}([0,1], \R_+)$  a probability
density function satisfying $g(z)= g(1-z)$ and $g(0)=0$. We define for a fixed $N\ge 2$  the probability measure
$$
\mu(du) := g\big(\{(N-1)u \}\big)du
$$ 
where $\{\cdot\}$ denotes the fractional part function. This defines an absolutely continuous  probability measure $\mu$ on the unit interval with an a.e. positive  continuous density. 

One checks that  $x^\star:=\big(\frac{k-1}{N-1} \big)_{k=1,\ldots,N}$ satisfies $\nabla L_N\big(x^\star\big)=0$ using the master equation derived from~\eqref{eq:grad2} and the obvious facts that $\Delta x_i = \frac{1}{N-1}$, $i=2,\ldots,N$, and $\int_0^1 \Psi_r(z)g(z) dz =  \int_0^1 \Psi_r(1-z)g(z) dz$.

Elementary computations starting from~\eqref{eq:Hess_i-1i} and~\eqref{eq:Hess_ii} show that the Hessian $\nabla^2L_{_{\!N}}(x^\star)$ is a symmetric  tridiagonal  $(N-2)\times (N-2)$ matrix  of the form 
\[
\nabla^2L_{_{\!N}}(x^\star)= (N-1)^{-(r-1)}A \quad \mbox{ with } \quad A=  \left[\begin{array}{cccccc}
2a&b&0&\cdots&\cdots&0\\
b&2a&b&0&\cdots&\vdots\\
0&\ddots&\ddots&\ddots&\ldots&0\\
\vdots&\ddots&b&\ddots&2a&b\\
0&\ldots&\ldots&0&b&2a
\end{array}\right] \!\in {\cal S}(N-2, \R). 
\]
with
$$
a= \int_0^1 \widetilde \vartheta_r (z)g(z)dz\quad \mbox{ and  } \quad  b= -  \int_0^1 \vartheta_r (z)g(z)dz <0
$$
(since both $\vartheta_r$ and $g$ are positive on $(0,1)$). It is classical background that (real) eigenvalues of such a symmetric tridiagonal matrix $A$ are 
\begin{equation}\label{eq:eigenvalues}
\lambda_k = 2\Big(a+b \cos \big(\tfrac{k\pi}{N-1}\big)\Big),\; k=1,\ldots, N-2,
\end{equation}
so that its lowest eigenvalue is $\lambda_{\min}(N)= 2\big(a+b\cos (\tfrac{\pi}{N-1})\big)$ (obtained with $k= 1$).

Using again that  $\vartheta_r -\widetilde \vartheta_r =-\Psi_r'$, 
\[
a +b= \int_0^1 (\widetilde \vartheta_r -\vartheta_r )(z)g(z)dz  = -\int_0^1 \Psi'_r(z) g(z)dz.
\]
Now, we note that $\Psi'_r(z) +\Psi'_r(1-z)= -\varpi''_r(z)$, $z\!\in [0,1]$, so that, taking advantage of the fact that $g(1-z)=g(z)$, we derive  
\[
\int_0^1 \Psi'_r(z) g(z)dz =  \int_0^{1} \frac{\Psi'_r(z) +\Psi'_r(1-z)}{2}g(z)dz= -\frac 12\int_0^{1}\varpi''_r(z) g(z)dz=-\int_0^{1/2}\varpi''_r(z) g(z) dz
\]
since both $g$ and $\varpi''_r$ are symmetric (w.r.t. $1/2$) on $[0,1]$. If we assumes  that $g$ is also differentiable on $(0,1)$, then  an integration by part  yields
\[
a+ b = \int_0^{1/2}\varpi''_r(z) g(z)dz =  g(1/2)\varpi'_r(1/2) -g(0)\varpi'(0)- \int_0^{1/2}\varpi'_r(z)g'(z)dz = - \int_0^{1/2}\varpi'_r(z)g'(z)dz
\]
since $g(0)=0$ and, by symmetry w.r.t $1/2$, $\varpi'_r(1/2) =0$. Finally setting 
\begin{equation}\label{eq:Ctrxple}
g= c_r \varpi_r
\end{equation}
so that $g$ is a probability density, one has 
$$
a+b=-c_r \int_0^{1/2}(\varpi'_r(z))^2dz<0,
$$
hence $\frac{a+b}{b}>0$ since $b<0$. Consequently, $\lambda_{\min}(N)<0$ for any $N$ large enough such that 
\[
\cos \big(\tfrac{\pi}{N-1}\big)>1- \frac{a+b}{b}.
\]
Then,   as some of the eigenvalues of the Hessian of $L_N$ at $x^\star$   are negative,  this point cannot be  a local minima of $L_N$. The function $L_N$ also has a local minima lying in ${\cal S}^{a,b}_{N}$ since ${\rm supp}(\mu)= [0,1]$, namely any $L^r$-optimal dual $N$-quantizer. Hence, uniqueness of the solution to the master equation  fails.

However, note that this only stands as a counter-example to uniqueness of solutions of  the master equation: the set of local minima may still be reduced to a single $N$-tuple. 

\bigskip
\noindent {\bf Remarks} $\bullet$ The conclusion of the above counterexample holds true for any symmetric probability density $g$ on $[0,1]$ such that $g(0)=0$ and 
\[
\int_0^{1/2}g'(z)\varpi'_r(z)dz>0.
\]
$\bullet$ If $\mu = U([a,b])$ then the master equation e.g. derived from~\eqref{eq:grad3} reads $(\Delta x_i)^{r} =  (\Delta x_{i+1})^{r} $, $i=2:N-1$ since $\int_{0}^{1}\Psi_r(z)dz= \int_{0}^{1}\Psi_r(1-z)dz$. Hence  $\Delta x_i = \frac{b-a}{N-1}$, $i=2:N$, so that one retrieves the fact that the unique $L^r$-optimal dual $N$-quantizer of $U([0,1])$ is always $x^{\star,N}=\big(\frac{i-1}{N-1}\big)_{i=1:N}$ for all $r\ge 1$. 

\smallskip
\noindent $\bullet$ If $r$ is an integer, then $\Psi_r(z) $ is  a polynomial function  with degree at most $r+1$. To be more precise one checks that its term of degree $r$ is always $0$ and that the coefficient $(-1)^r \big((-1)^r-1\big)z^{r+1}$ of its term of degree $r+1$ is $0$ if and only if  $r$ is even. Hence $d^0\, \Psi_r= r-1$ if $r$ is an even integer and $r+1$ if $r$ is an odd integer. 


\section{A Lloyd like algorithm for dual quantization in the quadratic case}\label{sec:llyod}
\paragraph{A fixed point formulation of the master equation.}

In this section, we specialize to the quadratic $r=2$ case and take advantage of this specialization to derive a more convenient expression of the distortion of the dual grid $\Gamma=\{a,x_2,\ldots,x_{N-1},b\}$ parametrized by $(x_2,\ldots,x_{N-1})\in {\cal S}_{_{\!N}}^{a,b}$ :
\begin{align}
\nonumber 
L_{_N}(x) & =\int_{[a,b]} \mu(d\xi)\int_0^1du |\xi-{\rm Proj}_{\Gamma}^{del}(\xi,u)\big|^2=  \int_{[a,b]} \mu(d\xi)\int_0^1du \big({\rm Proj}_{\Gamma}^{del}(\xi,u)\big)^2 - \int_\R \xi^2\mu(d\xi)\\
\nonumber
 &= \sum_{i=1}^{N-1}\int_{(x_i,x_{i+1}]}\mu(d\xi)\Big[  \tfrac{x_{i+1}-\xi}{x_{i+1}-x_i}x_i^2+\tfrac{\xi-x_{i}}{x_{i+1}-x_i} x_{i+1}^2\Big] - \int_\R \xi^2\mu(d\xi)\\
 \label{eq:dualweights0} &= \sum_{i=1}^{N-1}\int_{(x_i,x_{i+1}]}\hskip- 1cm \big((x_{i+1} +x_{i})\xi -x_{i+1} x_{i}\big)\mu(d\xi)-\int_\R \xi^2\mu(d\xi)\\
\label{eq:dualweights}
& =  \sum_{i=1}^{N-1} \left((x_i+x_{i+1})[K]_{x_{i}}^{x_{i+1}}-x_ix_{i+1}[F]_{x_{i}}^{x_{i+1}} \right)-\int_\R \xi^2\mu(d\xi),
\end{align}
where, we recall that $F(x)=\mu((-\infty,x])$ and $K(x)=\int_{(-\infty,x]}\xi\mu(d\xi)$ respectively denote the cumulative distribution function and the first partial moment of $\mu$ and for a function $g:\R \to \R$ and two real numbers  $x\le y$, $[g]_x^y= g(y)-g(x)$.
Then, it follows from~\eqref{eq:dualweights0} that
the mapping $ x\mapsto L_{_N}(x)$  is  continuously differentiable at $x$ when the distribution $\mu$ is atomless (i.e. $F$ and $K$ are continuous)
with
\begin{align}\label{eq:dualweights}
\partial_{x_i}  L_{_N}(x)=
[K]_{x_{i-1}}^{x_{i+1}}-\big(x_{i+1}[F]_{x_{i}}^{x_{i+1}}+x_{i-1}[F]_{x_{i-1}}^{x_{i}}\big),\quad i=2: N-1,
\end{align}

The master equation for optimal quadratic quantizers reads $x_1= a$, $x_{N}=b$ and 
\begin{equation}\label{eq:mastereq}
\nabla_{x_{2:N-1}}  L_{_N}(x)=0
\end{equation}
that is
\begin{equation}\label{eq:mastereq2}
[K]_{x_{i-1}}^{x_{i+1}} = x_{i+1}[F]_{x_{i}}^{x_{i+1}}+x_{i-1}[F]_{x_{i-1}}^{x_{i}},\quad i=2: N-1.
\end{equation}
%
%
%

\medskip Using Fubini's theorem for the second equality, we obtain that
\begin{align*}
  [K]_{x_{i-1}}^{x_{i+1}}&=\int_{(x_{i-1},x_{i+1}]}\int_{(x_{i-1},x_{i+1}]}\hskip-0.25cm \mbox{\bf 1}_{\{y<\xi\}}dy\,\mu(d\xi)+x_{i-1}[F]_{x_{i-1}}^{x_{i+1}}\\
  &=\int_{(x_{i-1},x_{i+1}]}\big(F(x_{i+1})-F(y)\big)dy+x_{i-1}[F]_{x_{i-1}}^{x_{i+1}}\\&=-\int_{x_{i-1}}^{x_{i+1}}F(y)dy+x_{i+1}F(x_{i+1})-x_{i-1}F(x_{i-1}).
\end{align*}
We deduce  the following more synthetic form for the master equation
\[
(x_{i+1}-x_{i-1})F(x_i) = \int_{x_{i-1}}^{x_{i+1}} F(\xi)d\xi, \quad i=2:N-1
\]
which may also be deduced from the case $r=2$ in \eqref{eq:grad3} using $\Psi'_2=1$ and performing a change of variables in each integral.
Equivalently, we have
\begin{align}
\label{eq:Tfixe2}
F(x_i) = \frac{\int_{x_{i-1}}^{x_{i+1}} F(\xi)d\xi}{x_{i+1}-x_{i-1}}, \quad i=2:N-1.
\end{align}
Assume from now on that the distribution $\mu$ is atomless with 
support $[a,b]$. Then $F:[a,b]\to [0,1]$ is   an increasing homeomorphism and  we may define its inverse $F^{-1}$ so  that the above  equation~\eqref{eq:Tfixe2}  can also be written as the fixed point equation  
\[
(x_2,\ldots,x_{N-1}) = T(x)= \big(T_2(x), \ldots,T_{{N-1}}(x)\big)
\]
where, for every $x \!\in  {\cal S}^{a,b}_{N}$,
\begin{equation}\label{eq:Lloydmap}
T(x) = \left( F^{-1}\left(   \frac{\int_{x_{i-1}}^{x_{i+1}} F(\xi)d\xi}{x_{i+1}-x_{i-1}}\right)\right)_{i=2:N-1}.
\end{equation} 
Since $F$ is continuous and increasing from $F(a)=0$ to $F(b)=1$, $T(x)\in{\cal S}^{a,b}_{N}$ for each $x\in{\cal S}^{a,b}_{N}$.
By Corollary~\ref{cor:unique}, any  quadratic optimal dual quantization grid at level $N$ is of the form $\{a,x_2,\ldots,x_{N-1},b\}$ with $(x_2,\ldots,x_{N-1})\in {\cal S}^{a,b}_{N}$ solution to the master equation $\nabla L_N(x)=0$ or equivalently fixed point of $T$. If $\mu$ admits a density $f$ with respect to the Lebesgue measure which is positive and $\log$-concave on the interval $(a,b)$ and vanishes outside, then there is a unique such point $
(x^\star_2,\ldots,x^\star_{N-1})$. 

One checks that $T$ {\em can be continuously extended to  the closure  $\overline  {\cal S}^{a,b}_{N}$ of  ${\cal S}^{a,b}_{N}$}. Indeed  $x_{i-1}=x_i <x_{i+1}$ or $x_{i-1}<x_i =x_{i+1}$, the extension of $T_i(x)$ is straightforward and    if $x_{i-1}=x_i =x_{i+1}$,  set $T_i(x)= x_i$ (in both cases with $x_1=a$ and $x_{N}=b$).

\medskip
From such a fixed point identity, one can devise an iterative fixed point procedure which can be seen as the counterpart of  so-called Lloyd's method~I procedure for dual quantization:
\begin{equation}\label{eq:DualLloyd}
x^{[\ell+1]}=T\big(x^{[\ell]}\big), \; \ell\ge 0, \quad x^{[0]}\!\in \overline{\cal S}^{a,b}_{N}.
\end{equation}

When $\mu$ admits a density $f$  positive and $\log$-concave on $(a,b)$ and vanishing outside, while proving that this procedure converges at a geometric rate to $x^\star$, we are going to check that $T$ admits a unique fixed point in ${\cal S}^{a,b}_{N}$ thus providing an alternative argument for the uniqueness statement in Theorem~\ref{thm:unique}.
 
\paragraph{Convergence of the dual Lloyd algorithm.} 

\medskip First note that if $F$ is  continuously differentiable on $(a,b)$ i.e. $\mu$ has a continuous density $f$ on $(a,b)$, then the mapping $T$ is itself continuously differentiable at any   $x\!\in{\cal S}^{a,b}_{N}$ with a Jacobian matrix $J_T(x) = \big[\frac{\partial T_i}{\partial x_j}(x)\big]_{2\le i,j\le N-1}$
where $\frac{\partial T_i}{\partial x_j}(x)= 0$ if $|i-j|\ge 2$ and 
\begin{equation}\label{eq:partialT1}
\frac{\partial T_i}{\partial x_{i-1}}(x) = \frac{1}{f\circ F^{-1}\Big(\frac{ \int_{x_{i-1}}^{x_{i+1}} F(\xi)d\xi}{x_{i+1}-x_{i-1}}\Big)}\times\frac{ \int_{x_{i-1}}^{x_{i+1}} (F(\xi)-F(x_{i-1}))d\xi }{(x_{i+1}-x_{i-1})^2} >0, \quad i= 3,\ldots,N-1,
\end{equation} 
and 
\begin{equation}\label{eq:partialT2}
\frac{\partial T_i}{\partial x_{i+1}}(x) = \frac{1}{f\circ F^{-1}\Big(\frac{ \int_{x_{i-1}}^{x_{i+1}} F(\xi)d\xi}{x_{i+1}-x_{i-1}}\Big)}\times \frac{ \int_{x_{i-1}}^{x_{i+1}} (F(x_{i+1})-F(\xi))d\xi }{(x_{i+1}-x_{i-1})^2}>0 ,\quad i= 2,\ldots,N-2.
\end{equation} 
The main result of this section is the following. 
\begin{Theorem}\label{thm:dualLloyd} Let $-\infty < a < b< +\infty$. Assume that $\mu$ admits a density $f$ with respect to the Lebesgue measure which is positive and $\log$-concave on the interval $(a,b)$ and vanishes outside. Then $T$ admits a unique fixed point $x^\star$ in $\overline {\cal S}^{a,b}_{N}$ and $x^\star\in{\cal S}^{a,b}_{N}$. Moreover, $\{a,x^\star_2,\ldots,x^\star_{N-1},b\}$ is the unique quadratic optimal dual quantization grid at level $N$ of $\mu$ and 
\begin{equation*}
   \exists \rho\in[0,1),\;\forall x^{[0]}\in\overline{\cal S}^{a,b}_{N},\forall \ell\in\N,\;|x^{[\ell]}-x^\star|_{\ell^\infty}\le \begin{cases}
     |x^{[0]}-x^\star|_{\ell^\infty}\rho^{\lfloor \ell/\widetilde N\rfloor}\mbox{ where $\widetilde N=\lceil\frac N2\rceil-1$}\\
     |x^{[0]}-x^\star|_{\ell^\infty}\rho^{\ell}\mbox{ if $f$ is strictly $\log$-concave on $(a,b)$.}
   \end{cases}
 \end{equation*}
               %
%
\end{Theorem}

As a preamble, we first establish an equivalence between two characterizations of {\em strong unimodality}, one being the $\log$-concavity of the density of  the distribution, whereas the other (see below) will be extensively used in what follows. Then we  establish a general result about fixed point of locally contracting transforms from a compact set in itself.

\begin{Lemma}\label{lem:logconcave}
   Let $-\infty < a < b< +\infty$ and let  $f: (a,b)\to (0,+\infty)$ be a positive probability density  on $(a,b)$  with cumulative distribution function $[a,b]\ni x\mapsto F(x)=\int_a^xf(\xi)d\xi$ and quantile function $F^{-1}$. The function $f$ is $\log$-concave (resp. strictly $\log$-concave) on $(a,b)$ iff $f\circ F^{-1}$ is concave (resp. strictly concave) on $(0,1)$.
\end{Lemma}
\noindent {\bf Proof.}
The cumulative distribution function $F$ being continuous and increasing on $[a,b]$ with $F(a)=0$ and $F(b)=1$, it admits a continuous and increasing inverse $F^{-1}:[0,1]\to[a,b]$.

Let us suppose that $f$ is $\log$-concave (resp. strictly $\log$-concave). Then $\log f$ is continuous and admits a non-increasing (resp. decreasing) right-hand derivative $(\log f)'_r$ as a real-valued concave (resp. strictly concave) function. By composition with the exponential, $f=\exp\circ\log f$ also admits a right-hand derivative equal to $f\times(\log f)'_r$. 
Since $f$ is continuous and positive, the function $F$ and its inverse $F^{-1}$ are continuously differentiable with respective derivatives $f$ and $\frac{1}{f\circ F^{-1}}$. We conclude that $f\circ F^{-1}$ admits a right-hand derivative equal to $\frac{f\times(\log f)'_r}{f}\circ F^{-1}=(\log f)'_r\circ F^{-1}$ which is non-increasing (resp. decreasing) as the composition of the non-increasing (resp. decreasing) function $(\log f)'_r$ with the increasing function $F^{-1}$. Therefore $f\circ F^{-1}$ is concave (resp. strictly concave).

When $f\circ F^{-1}$ is concave (resp. strictly concave), then this function is continuous and, by composition with the continuous function $F$, $f$ is continuous so that $F$ and $F^{-1}$ are continuously differentiable with respective derivatives $f$ and $\frac{1}{f\circ F^{-1}}$. Moreover $f=(f\circ F^{-1})\circ F$ admits a right-hand derivative equal to $f'_r=(f\circ F^{-1})'_r\circ F\times f$. Then $\log f$ admits a right-hand derivative equal to $\frac{f'_r}{f}=(f\circ F^{-1})'_r\circ F$ which is non-increasing (resp. decreasing) as the composition of the non-increasing (resp. decreasing) function $(f\circ F^{-1})'_r$ with the increasing function $F$. Therefore $\log f$ is concave (resp. strictly concave).
\hfill $\Box$

\begin{Proposition}\label{prop:locally-contract} Let $K$ be a convex compact subset of $\R^d$, equipped with a norm $\|\cdot \|$,  and let $T: K\to K$ be a $\|\cdot \|-1$-Lipschitz continuous mapping such that for some $k\in\N^*$, the mapping $T^k$ obtained by iterating $T$ $k$-times satisfies
  $$\rho_k:=\sup_{y \in K,\,y\neq y^\star} \tfrac{\|T^k(y)-y^\star\|}{\|y-y^\star\|}<1$$ for some fixed point $y^\star$ of $T$ (the set of fixed points is non empty by Brouwer's theorem). Then $y^\star$ is the unique fixed point of $T$ and for every $y_0\!\in K$, the sequence recursively defined for $n\in\N$ by 
  $y_{n+1}= T(y_n)$ geometrically converges to $y^\star$ :
$$\forall n\in\N,\;\|y_n-y^\star\|\le \|y_0-y^\star\|\rho_k^{\lfloor n/k\rfloor}.$$
\end{Proposition}

\noindent {\bf Proof.}
The inequality $\|T^k(y)-y^\star\|\le \rho_k\|y-y^\star\|$ valid for each $y\in K$ with $\rho_k<1$ implies that every fixed point of $T^k$ and therefore of $T$ is equal to $y^\star$. Moreover, using that $y^\star$ is a fixed point of $T$ then the $1$-Lipschitz property of $T$ and last the previous inequality, we obtain that
\begin{align*}
   \|y_n-y^\star\|&=\|T^{n-\lfloor n/k\rfloor k}(y_{\lfloor n/k\rfloor k})-T^{n-\lfloor n/k\rfloor k}(y^\star)\|\le \|T^k(y_{(\lfloor n/k\rfloor-1) k})-y^\star\|\\&\le \rho_k\|T^k(y_{(\lfloor n/k\rfloor-2) k})-y^\star\|\le \rho_k^{\lfloor n/k\rfloor}\|y_0-y^\star\|.
\end{align*}
\hfill$\Box$


 \noindent{\bf Proof of Theorem~\ref{thm:dualLloyd}.} Since $K=\overline {\cal S}_{_{\!N}}^{a,b}$ is a convex and compact subset of $\R^{N-2}$, by Brouwer's fixed-point theorem, the set of fixed points of the continuous map $T:\overline {\cal S}_{_{\!N}}^{a,b}\to \overline {\cal S}_{_{\!N}}^{a,b}$ is non empty. Let $x\in\overline {\cal S}_{_{\!N}}^{a,b}$ and $x_1=a$, $x_N=b$. If $T_i(x)=x_{i-1}$ or $T_{i}(x)=x_{i+1}$ for some $i=2,\ldots,N-1$, then, since $F$ is increasing and continuous on $[a,b]$, $x_{i-1}=x_i=x_{i+1}$. If moreover $T(x)=x$, then one deduces that $T_{i+1}(x)=x_i$ if $i\le N-2$ and $T_{i-1}(x)=x_i$ if $i\ge 3$, so that by induction $x_1=x_2=\ldots=x_{N-1}=x_{N}$ which contradicts $x_1=a<b=x_N$. Hence the non-empty set of fixed points of $T$ is included in  ${\cal S}_{_{\!N}}^{a,b}$.
Let $x^\star\in{\cal S}_{_{\!N}}^{a,b}$ be one of these fixed points. We are going to check that the assumptions of Proposition~\ref{prop:locally-contract} are satisfied with $y^\star=x^\star$ and with $k=1$ in the strictly $\log$-concave case and $k=\widetilde N$ is the $\log$-concave case. The conclusions but the link between $x^\star$ and the unique quadratic optimal dual quantization grid at level $N$ of $\mu$ then follow from this proposition. This link is a consequence of Corollary~\ref{cor:unique} and the fact that  $x\in{\cal S}_{_{\!N}}^{a,b}$ is a critical point of $L_N$ iff it is a fixed point of $T$. 
 

\smallskip
\noindent -- {\em Srictly $\log$-concave setting}. It follows from~\eqref{eq:partialT1} and~\eqref{eq:partialT2} that, for every $x\!\in {\cal S}_{_{\!N}}^{a,b}$ and $i=3,\ldots,N-2$,
\[
\frac{\partial T_i}{\partial x_{i-1}}(x)+\frac{\partial T_i}{\partial x_{i+1}}(x)=  \frac{1}{f\circ F^{-1}\Big(\frac{ \int_{x_{i-1}}^{x_{i+1}} F(\xi)d\xi}{x_{i+1}-x_{i-1}}\Big)}\times\frac{F(x_{i+1})-F(x_{i-1}) }{x_{i+1}-x_{i-1}}.
\]
As $f$ is strictly $\log$-concave, $f\circ F^{-1}$  is strictly concave by Lemma~\ref{lem:logconcave}, hence
it follows from Jensen's inequality that, under the convention $x_1=a$ and $x_N=b$, for every $i=2,\ldots,N-1$
\begin{equation}
   f\circ F^{-1}\left(\frac{ \int_{x_{i-1}}^{x_{i+1}} F(\xi)d\xi}{x_{i+1}-x_{i-1}}\right) >\frac{ \int_{x_{i-1}}^{x_{i+1}} f(\xi)d\xi}{x_{i+1}-x_{i-1}} = \frac{F(x_{i+1})-F(x_{i-1}) }{x_{i+1}-x_{i-1}}\label{eq:strictconcav}
\end{equation}
with a strict inequality since the probability measure $\mbox{\bf 1}_{[x_{i-1},x_{i+1}]}(\xi)d\xi$ is not a Dirac mass and $F$ is not constant over $[x_{i-1},x_{i+1}]$.  As a consequence, for every $x\!\in {\cal S}^{a,b}_{N}$,  
\begin{align}\label{eq:SumJ_T}
  &\forall i=3,\ldots,N-2,\;\frac{\partial T_i}{\partial x_{i-1}}(x)+\frac{\partial T_i}{\partial x_{i +1}}(x) <1,\\
  &\frac{\partial T_2}{\partial x_{3}}(x)=\frac{\int_a^{x_3}(F(x_3)-F(\xi))d\xi }{(x_{3}-{a})^2f\circ F^{-1}\Big(\frac{ \int_{a}^{x_{3}} F(\xi)d\xi}{x_{3}-a}\Big)}<\frac{1}{f\circ F^{-1}\Big(\frac{ \int_{a}^{x_{3}} F(\xi)d\xi}{x_{3}-a}\Big)}\times\frac{F(x_3)-F(a) }{x_{3}-a}<1,\label{eq:SumJ_T2}\\
  &\frac{\partial T_{N-1}}{\partial x_{N-2}}(x)=\frac{\int^b_{x_{N-2}}(F(\xi)-F(x_{N-2}))d\xi }{(b-x_{N-2})^2f\circ F^{-1}\Big(\frac{ \int^b_{x_{N-2}}F(\xi)d\xi}{b-x_{N-2}}\Big)}<\frac{1}{f\circ F^{-1}\Big(\frac{ \int^b_{x_{N-2}}F(\xi)d\xi}{b-x_{N-2}}\Big)}\times\frac{F(b)-F(x_{N-2}) }{b-x_{N-2}}<1.\label{eq:SumJ_TN-1}\end{align}
Since $\frac{\partial T_i}{\partial x_j}(x)\ge 0$  with equality if $|i-j|\ge 2$, it is easy to deduce that, if $\R^{N-2}$ is equipped with the $\ell^{\infty}$-norm $|\cdot|_{\ell^\infty}$,
\begin{equation}\label{eq:NormJ_T}
\forall\, x\! \in {\cal S}^{a,b}_{N},\quad \vertiii{J_T(x)}_{\ell^\infty}   <1
\end{equation}
where $\vertiii{\cdot}_{\ell^\infty} $ denotes the operator norm with respect to the $\ell^\infty$-norm.  

Now let $x\!\in \overline{\cal S}^{a,b}_{N}$ and $y\!\in{\cal S}^{a,b}_{N}$ with $y\neq x$. Then for each $t\in[0,1)$, $tx+(1-t)y\in{\cal S}^{a,b}_{N}$ and $t\mapsto T(tx+(1-t)y)$ is continuous on $[0,1]$ and differentiable on $[0,1)$ so that, with the integrability consequence of \eqref{eq:NormJ_T}, 
\begin{equation}
   T(x)-T(y) =\int_0^1J_T\big(t x+(1-t)y\big) (x-y)dt,\;\forall (x,y)\!\in \overline{\cal S}^{a,b}_{N}\times{\cal S}^{a,b}_{N}.\label{integdif}
\end{equation}
By the triangle inequality for integrals,  the definition of the $\vertiii{\cdot}_{\ell^{\infty}}$-norm and~\eqref{eq:NormJ_T}, one deduces that
\begin{equation}\label{eq:Tquasi-contract}
\big| T(x)-T(y) \big|_{\ell^\infty} \le \int_{0}^1 \vertiii{J_T\big(t x+(1-t)y\big)}_{\ell^\infty} dt \,|x-y|_{\ell^\infty}<|x-y|_{\ell^\infty}.
\end{equation}
Approximating $y\in\overline{\cal S}^{a,b}_{N}$ by a sequence of elements in ${\cal S}^{a,b}_{N}$, we deduce that $T$ is $|.|_{\ell^\infty}$ $1$-Lipschitz continuous. 
Moreover, for the choice $y$ equal to the fixed point $x^\star$ of $T$, \eqref{integdif} writes
\begin{equation}
\forall x\in\overline{\cal S}^{a,b}_{N},\quad T\big(x\big)-x^\star= A_{x}^{[1]} \big(x-x^\star\big)\quad \mbox{with} \quad   A^{[1]}_{x}:= \int_0^1 J_T\big(x^\star+t(x-x^\star)\big)dt.\label{diftxx*}
\end{equation}
Since $x\mapsto A^{[1]}_{x}$ is continuous on the compact set $\overline{\cal S}^{a,b}_{N}$ and $$\forall x\in \overline{\cal S}^{a,b}_{N},\;\vertiii{A^{[1]}_{x}}_{\ell^\infty}\le \int_0^1\vertiii{J_T\big(t x+(1-t)x^\star\big)}_{\ell^\infty} dt<1,$$ we have $\sup_{x\in \overline{\cal S}^{a,b}_{N}}\vertiii{A^{[1]}_{x}}_{\ell^\infty}<1$. With \eqref{diftxx*}, we deduce that the hypotheses of Proposition~\ref{prop:locally-contract} are satisfied with $k=1$ and $\rho_1=\sup_{x\in \overline{\cal S}^{a,b}_{N}}\vertiii{A^{[1]}_{x}}_{\ell^\infty}$.

\smallskip
\noindent -- {\em $\log$-concave setting}.

The main difference is that the inequality in \eqref{eq:strictconcav} is no longer strict. As a consequence, in \eqref{eq:SumJ_T}-\eqref{eq:SumJ_TN-1}, $<1$ should now be replaced by $\le 1$ so that  $\vertiii{J_T(x)}_{\ell^{\infty}}\le 1$. This still ensures that $T$ is $|{\cdot}|_{\ell^{\infty}}$-$1$-Lipschitz continuous on $\overline{\cal S}_{_{\!N}}^{a,b}$ and \eqref{diftxx*} still holds. To overcome the lack of strict contraction of $\R^{N-2}\ni u\mapsto J_T(x)u$, we are going to take advantage of the inequalities $\frac{\partial T_2}{\partial x_{3}}<1$ and $\frac{\partial T_{N-1}}{\partial x_{N-2}}<1$ still valid on ${\cal S}_{_{\!N}}^{a,b}$ since the first inequalities in 
\eqref{eq:SumJ_T2}-\eqref{eq:SumJ_TN-1} remain strict.

For $k\ge 1$ we can iterate \eqref{diftxx*} to obtain
\begin{equation}
 \forall x\in\overline{\cal S}^{a,b}_{N},\;T^k\big(x\big)-x^\star= A_{x}^{[k]}A_{x}^{[k-1]}\ldots A_{x}^{[1]}\big(x-x^\star\big)\label{diftxx*it}  
\end{equation}
where, for $\ell\ge 1$, $A_{x}^{[\ell]}= \int_0^1 J_T\big(x^\star+t(T^{l-1}(x)-x^\star)\big)dt$
is a tridiagonal matrix  $\big[a^{[\ell]}_{ij}\big]_{2\le i,j\le N-1}$ satisfying
\begin{equation}
\label{eq:A1}
 \begin{array}{lll} &a^{[\ell]}_{ii} = 0, \;\; i=2,\ldots,N-1,&\quad  a^{[\ell]}_{i i\pm1} >0, \, a^{[\ell]}_{ii-1}+a^{[\ell]}_{ii+1}\le 1, \;\; i=3,\ldots,N-2, \\
\mbox{and } &0<a^{[\ell]}_{23}, \;a^{[\ell]}_{N-1N-2} <1.
 \end{array} 
\end{equation}
Let $u= (u_2,\ldots,u_{_{N-2}})\!\in \R^{N-2}$ with  $|u|_{\ell^\infty} >0$. Let $I_1= \big\{i: |A^{[1]}_xu|=|u|_{\ell^\infty} \big\}$. It is clear that $2,N-1\!\notin I_1$ since $\big|(A^{[1]}_xu)_2 \big|= a^{[1]}_{23}|u_3| <|u_3|\le |u|_{\ell^\infty} $ (idem for the other term).  Now  let $I_2=    \big\{i: \big|(A^{[2]}_xA^{[1]}_xu)_i \big|=|u|_{\ell^\infty}\big\}$. It is still clear that $2, N-1\!\notin I_2$. Now 
\[
\big(A_x^{[2]}A^{[1]}_xu\big)_3  = a^{[2]}_{32} \big(A^{[1]}_xu\big)_2  + a^{[2]}_{34} \big(A_x^{[1]}u\big)_4 .
\]
As $a^{[2]}_{32},\, a^{[2]}_{34} >0$ and $a^{[2]}_{32}+a^{[2]}_{34} \le 1$ and $2\notin I_1$, $|(A_x^{[2]}A_x^{[1]}u)_3  |<|u|_{\ell^\infty}$. One shows likewise that $N-2\!\notin I_2$. Then one shows  the same way round by induction that,
$2,\ldots, k+1$, $N-(k+1),\ldots, N-1\!\notin I_k= \big\{i: \big|(A_x^{[k]}\cdots A_x^{[1]}u)_i \big|= |u|_{\ell^\infty} \big\}$ which implies that $I_{\lceil \frac{N}{2}\rceil-1}= \varnothing$ i.e.
$\big | A_x^{[\lceil \frac{N}{2}\rceil-1]}\cdots A_x^{[1]}u\big|_{\ell^{\infty}}<|u|_{\ell^\infty}$. Consequently setting $\widetilde N =\lceil \frac{N}{2}\rceil -1$, we have  \begin{equation*}
   \forall x\in \overline{\cal S}_{_{\!N}}^{a,b},\;\vertiii{A^{[\widetilde N]}_x\cdots A^{[1]}_x}_{\ell^{\infty}}<1.
 \end{equation*}
Note that a more quantitative bound in terms of the coefficients of the matrices is derived in Proposition~\ref{prop:quantbound} below.
 With the continuity of $x\mapsto A^{[\widetilde N]}_x\cdots A^{[1]}_x$ over the compact set $\overline{\cal S}_{_{\!N}}^{a,b}$ and \eqref{diftxx*it}, we deduce that the hypotheses of Proposition~\ref{prop:locally-contract} are satisfied with $k=\widetilde N$ and $\rho_{\widetilde N}=\sup_{x\in \overline{\cal S}_{_{\!N}}^{a,b}}\vertiii{A^{[\widetilde N]}_x\cdots A^{[1]}_x}_{\ell^{\infty}}<1$.  \hfill $\Box$
 
To evaluate in sharper way $\rho_{\widetilde N}$, one might rely on the following Proposition which provides a  quantitative bound for the $\vertiii{\cdot}_{\ell^{\infty}}$-norm of a product of our  tridiagonal matrices of interest.
 
\begin{Proposition}[Quantitative bound]\label{prop:quantbound}
   Let $N\ge 3$ and set $\tilde N=\lceil \frac N2\rceil -1$.     Let  $A^{[\ell]}=[a^{[\ell]}_{ij}]_{2\le i,j\le N-1}$,  $l\in\{1,\hdots,\tilde N \}$,  be tridiagonal matrices whose entries satisfy the above conditions~\eqref{eq:A1}.
Then \begin{align*}
   \vertiii{A^{[\tilde N]}\cdots A^{[1]}}_{\ell^\infty}\le &\max_{2\le i\le\tilde N}\bigg(1-a_{ii-1}^{[\tilde N]}a_{i-1i-2}^{[\tilde N-1]}\cdots a_{32}^{[\tilde N+3-i]}(1-a^{[\tilde N+2-i]}_{23})\bigg)\\&\vee\max_{\lfloor\frac{N}{2}\rfloor\le i\le N-1}\bigg(1-a^{[\tilde N]}_{ii+1}a^{[\tilde N-1]}_{ii+2}\cdots a^{[\tilde N+i+2-N]}_{N-2N-1}(1-a^{[\tilde N+1+i-N]}_{N-1N-2})\bigg)<1.
\end{align*}
\end{Proposition}
\noindent{\bf Proof.}
For a matrix $B=[b_{ij}]_{1\le i\le n,\;1\le j\le d}\in\R_+^{n\times d}$ and for $y,z\in\R^d$ such that $|y_i|\le z_i$, $i=1:d$, we have
$$\|By\|_{\ell^\infty}=\max_{1\le i\le n}\bigg|\sum_{j=1}^d b_{ij}y_j\bigg|\le \max_{1\le i\le n}\sum_{j=1}^d b_{ij}z_j=\|Bz\|_{\ell^\infty}.$$
Therefore denoting by ${\mathbf 1}$ the vector in $\R^{N-2}$ with all coordinates equal to $1$, we have $\vertiii{A^{[\tilde N]}\cdots A^{[1]}}_{\ell^\infty}\le\|A^{[\tilde N ]}\cdots A^{[1]}{\mathbf 1}\|_{\ell^\infty}$. To conclude, we check by induction on $k\in\{1,\cdots,\tilde N \}$ that the entry $(A^{[k]}\cdots A^{[1]}{\mathbf 1})_i$ is nonnegative and not greater than
\begin{equation*}
   \begin{cases}a^{[k]}_{ii+1}+a^{[k]}_{ii-1}\Big(a^{[k-1]}_{i-1i}+a^{[k-1]}_{i-1i-2}\big(a^{[k-2]}_{i-2i-1}+a^{[k-2]}_{i-2i-3}(\hdots+a^{[k+3-i]}_{32}(a^{[k+2-i]}_{23}+0))\big)\Big)\mbox{ if } 2\le i\le k+1\\
      \;1\hskip 11,5cm \mbox{ if }\quad k+2\le i\le N-(k+1)\\
      a^{[k]}_{ii-1}+a^{[k]}_{ii+1}\Big(a^{[k-1]}_{i+1i}+a^{[k-1]}_{i+1i+2}\big(a^{[k-2]}_{i+2i+1}+a^{[k-2]}_{i+2i+3}(\hdots+a^{[k+i+2-N]}_{N-2N-1}(a^{[k+i+1-N]}_{N-1N-2}+0))\big)\Big)  \mbox{ if } N-k\le i\le N-1
   \end{cases}
 \end{equation*}
 where, by the assumptions made on the entries of the matrices, by induction on $i\in\{2,\cdots,k+1\}$,
\begin{align*}
   a^{[k]}_{ii+1}+a^{[k]}_{ii-1}\Big(a^{[k-1]}_{i-1i}+a^{[k-1]}_{i-1i-2}\big(a^{[k-2]}_{i-2i-1}+a^{[k-2]}_{i-2i-3}&(\hdots+a^{[k+3-i]}_{32}(a^{[k+2-i]}_{23}+0))\big)\Big)\\&\le 1-a^{[k]}_{ii-1}a^{[k-1]}_{i-1i-2}\cdots a^{[k+3-i]}_{32}\big(1-a^{[k+2-i]}_{23}\big)<1
\end{align*} 
and, by backward induction on $i\in\{N-k,\cdots,N-1\}$,
\begin{align*}
   a^{[k]}_{ii-1}+a^{[k]}_{ii+1}\Big(a^{[k-1]}_{i+1i}+a^{[k-1]}_{i+1i+2}\big(a^{[k-2]}_{i+2i+1}+a^{[k-2]}_{i+2i+3}&(\hdots+a^{[k+i+2-N]}_{N-2N-1}(a^{[k+i+1-N]}_{N-1N-2}+0))\big)\Big)\\
   &\le 1-a^{[k]}_{ii+1}a^{[k-1]}_{i+1i+2}\cdots a^{[k+i+2-N]}_{N-2N-1}\big(1-a^{[k+i+1-N]}_{N-1N-2}\big)<1.
\end{align*}
In the induction on $k$, we use the first bound to get that $(A^{[k+1]}\cdots A^{[1]}{\mathbf 1})_i\le a^{[k+1]}_{ii+1}+a^{[k+1]}_{ii-1}(A^{[k]}\cdots A^{[1]}{\mathbf 1})_{i-1}$ for $3\le i\le k+1$ and the second one to get that $(A^{[k+1]}\cdots A^{[1]}{\mathbf 1})_i\le a^{[k+1]}_{ii-1}+a^{[k+1]}_{ii+1}(A^{[k]}\cdots A^{[1]}{\mathbf 1})_{i+1}$ for $N-k\le i\le N-2$.
\hfill$\Box$

\medskip

\bigskip
\noindent {\bf Remarks.} $\bullet$ Note that, as all the entries of  the matrix $A^{[1]}_{x}$ such that $T(x)-x^\star= A^{[1]}_x(x-x^\star)$ are non negative, it is clear that if $x\ge x^\star$ (resp. $x\le x^\star$) {\em componentwise }then $T(x) \ge x^\star$ (resp. $ T(x)\le x^\star)$ {\em componentwise} so that if $x^{[0]}\ge x^\star$ (resp. $x^{[0}\le x^\star$) {\em componentwise} then the whole sequence $x^{[\ell]}$ will satisfy the same inequality.

%

\smallskip
\noindent  $\bullet$ This theorem shows the geometric convergence of this dual Lloyd   procedure for the ($\log$-concave) truncated exponential distributions toward its unique {\em quadratic} optimal dual quantizer. A specific family of procedures which works for the search of the  $L^r$-optimal dual quantizer of  power distributions for any $r\ge 1$ and  truncated exponential  distributions for $r\in\{1,2\}$  is  developed in the next section. 

\smallskip
\noindent  $\bullet$ 
In the $\log$-concave example of the uniform distribution, say on the unit interval $[0,1]$, one checks that the mapping $T$ is affine and reads on $\overline {\cal S}_{_{\!N}}^{a,b}$
\[
T_i (x) =  \frac{x_{i-1}+x_{i+1}}{2}, \quad i=2,\ldots,N-1,
\]
with the convention $x_1=0$ and $x_{N}= 1$ so that $T(x)=Ax +b$ with $A=[a_{ij}]_{2\le i,j\le N-1}$ satisfying $a_{i, i\pm 1}= \frac 12$, $i=3,\ldots,N-2$, $a_{2,3}=\frac 12=a_{N-1,N-2}$ and $a_{ij}= 0$ otherwise and $b= \frac 12 (0, \cdots, 0,1)^*$ (with $N-2$ components). The eigenvalues of the symmetric matrix $A$ are $\big(\cos(\frac{k\pi}{N-1})\big)_{1\le k\le N-2}$ (see~\eqref{eq:eigenvalues}). Therefore, for the Euclidean norm, $$\forall x\in\R^{N-2},\;|Ax|\le \cos\left(\frac{\pi}{N-1}\right)
|x|. $$
Hence the sequence $x^{[\ell+1]}= T\big(x^{[\ell]}\big)$, $\ell  \ge 0$, converges toward the unique equilibrium point $\big(  \frac{k-1}{N-1} \big)_{k=2,\ldots, N-1}$ with the geometric rate $\cos\left(\frac{\pi}{N-1}\right)$ uniformly with respect to  $x^{[0]}\in\overline{{\cal S}}^{a,b}_{N}$.

\section{Computation of one-dimensional dual grids for specific distributions}\label{sec:specific}
\subsection{Power distributions on compact intervals}Let $\mu$ be a probability measure compactly supported on $[0,1]$ with density $f$ and such that $0<\mu([0,x])<1$ for all $x\in(0,1)$. Then $x^{(N)}_1=0$ and $x^{(N)}_N=1$. To characterize the other points $(x_i^{(N)})_{2\le i\le N-1}$ in the optimal dual grid when $N\ge 3$, we use the master equations written with the expression~\eqref{eq:grad2} of the gradient of the distortion :
\begin{align*}
(\Delta x^{(N)}_{i+1})^{r}  \int_0^{1}\Psi_r(1-z) f(x_{i}^{(N)} +z\Delta x^{(N)}_{i+1}\big) dz=(\Delta x^{(N)}_i)^{r} \int_{0}^{1}\Psi_r(z) f\big(x_{i-1}^{(N)} +z\Delta x^{(N)}_i\big) dz,\;i=2:N-1
\end{align*}
where $\Psi_r(z)=(r-1)z(1-z)^r+z^{r+1}+r z^2(1-z)^{r-1}$.
For power distributions, $f(x)=\alpha x^{\alpha-1}$ for some $\alpha>0$, so that dividing the master equation by $\alpha (x_i^{(N)})^{r+\alpha-1}$ yields

\begin{align}
\label{master2}  \left(\frac{\Delta x^{(N)}_{i+1}}{x_{i}^{(N)}}\right)^{r}\int_0^1\Psi_r(1-z)&\left(1+\frac{z\Delta x^{(N)}_{i+1}}{x^{(N)}_i}\right)^{\alpha-1} dz\\
\nonumber &=\left(\frac{\Delta x^{(N)}_{i}}{x_{i}^{(N)}}\right)^{r}\int_0^1\Psi_r(z)\left(1+\frac{(z-1)\Delta x^{(N)}_i}{x^{(N)}_i}\right)^{\alpha-1} dz,
  \;  i=2:N-1.
\end{align}
We are going to check that the ratios $\lambda_i=\frac{x^{(N)}_{i}}{x^{(N)}_{i+1}}$ do not depend on $N\ge i+1$ (they of course depend on $r\ge 1$ but do not make this dependence explicit in the notation). This is a consequence of the equality $\frac{x^{(N)}_1}{x^{(N)}_2}=0$ valid for each $N\ge 2$ and which yields $\lambda_1=0$. Since $x^{(N)}_N=1$, we then have $x^{(N)}_{i}=\prod_{j=i}^{N-1}\lambda_{j}$ for $i=1:N-1$ and even for $i=N$ under the usual convention $\prod_{j=N}^{N-1}\lambda_j=1$. Performing the change of variable $y=1-z$ in the integral in the right-hand side of~\eqref{master2}, we obtain
\begin{equation}
\chi_r(\lambda_i^{-1}-1)=\chi_r(\lambda_{i-1}-1)\label{iterlam}   
\end{equation} where
\begin{equation*}
   \chi_r(x)=\begin{cases}
     x^{r}\int_0^1\Psi_r(1-z)(1+z x)^{\alpha-1}dz\;\mbox{ if }\;x\ge 0\\
     (-x)^{r}\int_0^1\Psi_r(1-z)(1+z x)^{\alpha-1}dz\;\mbox{ if }\;x\!\in[-1,0].
   \end{cases}
 \end{equation*}
 To conclude that starting from $\lambda_1=0$, the values of $\lambda_i$ can be computed inductively for $i=2:N-1$ from this equation, it is enough to check that $\chi_r$ is one to one on the interval $(0,+\infty)$ where $x^{(N)}_{i+1}/x^{(N)}_i-1$ stands. For $x\in(0,+\infty)$, we have
 $$
 \chi'_r(x)=x^{r-1}\int_0^1\Psi_r(1-z)(r+(r+\alpha-1)z x)(1+z x)^{\alpha-2}dz
 $$
 where the right-hand side is positive since $r\ge 1$ and $\Psi_r$ is non-negative. Notice that we obtain uniqueness for the master equation and therefore uniqueness of the $L^r$-optimal dual quantization grid at level $N$ even when $\alpha\in(0,1)$ and the density is not $\log$-concave.  In the quadratic $r=2$ case, since $\Psi_2(z)=z$, $\alpha\chi_2(x)=\frac{(1+x)^{1+\alpha}}{1+\alpha}-x-\frac{1}{1+\alpha}$.
 Of course when for $a<b$, $\mu$ admits the density $\mbox{\bf 1}_{[a,b]}(x)\frac{\alpha(x-a)^{\alpha-1}}{(b-a)^\alpha}$ (resp. $\mbox{\bf 1}_{[a,b]}(x)\frac{\alpha(b-x)^{\alpha-1}}{(b-a)^\alpha}$) then for $i=1:N$, $x_i^{(N)}=a+(b-a)\prod_{j=i}^{N-1}\lambda_j$ (resp. $x_{N+1-i}^{(N)}=b-(b-a)\prod_{j=i}^{N-1}\lambda_j$).

\medskip
 \noindent {\bf Numerical example.} The optimal quadratic dual $10$-quantizer of $\mu(dx)\!=\!{\mathbf 1}_{[0,1]}(x)\frac{dx}{2\sqrt{x}}$ is given by 
 $$
 \{0,0.0744614,       0.1675381,      0.2704687,      0.3804786,      0.4961058,      0.6164311,      0.7408177,      0.868795,   1\}.
 $$ Notice that even if the density is $\log$-convex on the interval $(0,1)$, the Llyod-like iterative algorithm introduced in Section~\ref{sec:llyod} still numerically converges to the corresponding unique solution to the master equation in ${\cal S}^{0,1}_{10}$.

The derivation of an equation not depending on $x_i^{(N)}$ relating $\lambda_i$ to $\lambda_{i-1}$ was permitted by the key structure condition
\begin{equation}
   \forall (x,y)\in(0,1]\times[0,1],\;f(y)=g(x)h\left(\frac{y}{x}\right)\label{structcfd}
 \end{equation}
 satisfied by the power density for $g(x)=\alpha x^{\alpha-1}$ and $h(z)=z^{\alpha-1}$. Under this structure condition,~\eqref{iterlam} remains valid for the following generalized definition of $\chi_r$ :
\begin{equation*}
   \chi_r(x)=\begin{cases}
     x^{r}\int_0^1\Psi_r(1-z)h(1+z x)dz\mbox{ if }x\ge 0\\
     (-x)^{r}\int_0^1\Psi_r(1-z)h(1+z x)dz\mbox{ if }x\in[-1,0].
   \end{cases}
 \end{equation*}
 When the functions $g$ and $h$ are differentiable, the structure condition is only satisfied by   power distributions. Indeed, differentiating with respect to $x$ in~\eqref{structcfd}, we obtain
 $\frac{x g'(x)}{g(x)}=\frac{\frac{y}{x}h'(\frac{y}{x})}{h(\frac{y}{x})}$. We deduce that the two functions $\frac{x g'(x)}{g(x)}$ and $\frac{z h'(z)}{h(z)}$ are both equal to some constant $\alpha-1$. Then $g(x)\propto x^{\alpha-1}$ and $h(z)\propto z^{\alpha-1}$ so that $f(y)\propto y^{\alpha-1}$.

\smallskip
We could also assume that $f(y)=g(x)h(y-x)$, which is typically satisfied when $f(x)=\frac{|\lambda| e^{\lambda x}}{e^\lambda -1}$ for $\lambda\in\R^*$ (we may then choose $h(z)=e^{\lambda z}$) but then it is not so easy to decouple the use of the two boundary conditions $x^{(N)}_1=0$ and $x^{(N)}_N=1$ which permits an inductive resolution of the master equation under the former structure condition~\eqref{structcfd}. We are nevertheless able to design an almost explicit procedure for these truncated exponential distributions, at least when $\lambda>0$ and $r\in\{1,2\}$.

\subsection{Truncated exponential distributions}  Let $\mu(dx)={\mathbf 1}_{[a,b]}(x)\frac{\lambda e^{-\lambda (x-a)}}{1-e^{-\lambda b}}dx$ be a truncated exponential distribution with parameter $\lambda >0$ on $[a,b]$, $\infty  <a < b<+\infty $. Note that if $\lambda <0$, it suffices to solve the problem for  the image $\tilde \mu ={\mathbf 1}_{[-b,-a]}(x)\frac{|\lambda| e^{-|\lambda| (x+a)}}{e^{|\lambda| b}-1}dx$ of $\mu$ by the linear transform $x\mapsto -x$ and transport the resulting dual quantizer by this involution.   

The distribution $\mu$ is a $\log$-concave distribution (though not strictly $\log$-concave) so that, for every $r\ge1$,  the    $L^r$-optimal dual quantizer, solution to the $L^r$-master equation~\eqref{eq:grad2}  is unique at every level $N\ge 3$. Let $N\ge 3$ and let $x=(x_1,\ldots x_{N})$, $x_1 =a$, $x_{N}=b$ and $\Delta x_i = x_i-x_{i-1}$, $i=2,\ldots, N$. 
The master equation~\eqref{eq:grad2} reads  
\begin{equation}\label{eq:expeq}
\Phi_r(\lambda \Delta x_{i+1}) =\Phi_r(-\lambda \Delta x_i), \; i=2,\ldots,N-1, \;x_1=a, \; x_{N}=b,
\end{equation}
with 
\[
\Phi_r(x) =|x|^re^{- x}\int_0^1 \Psi_r(z) e^{x z} dz
, \quad  x\!\in \R.
\]
If $x^{N,\lambda,a,b}$ denotes the   solution to this equation (where the dependence on $r\ge 1$ is not made explicit), one easily checks, taking advantage of uniqueness of the solution, that
\[
x^{N, \lambda,a, b} =a\cdot\mbox{\bf 1}+ \frac{1}{\lambda}x^{N,1,0,\lambda (b-a)}
\]
so we only need to solve the equation with $\lambda =1$ and $a=0$ with limit condition $x^{N, 1,0, b}_{N} =b$.

\smallskip
-- {\em Quadratic case} ($r=2$). We first consider the quadratic case $r=2$, most commonly (sic) used in applications. Then, $\Psi_2(z) = z$, so that
\[
\Phi_2(x) =x^2e^{- x}\int_0^1 z e^{x z} dz = e^{-x}-1+x, \quad  x\!\in \R.
\]
 
The function $\Phi_2 $ is $C^1$ and if we set $\check \Phi_2 (x) = \Phi_2(-x)$, then  $\Phi_{2|\R_+}$ and $\check \Phi_{2|\R_+}$  are both    increasing $C^1$-diffeomorphisms of $\R_{+}$ and (with an obvious abuse of notation) Equation~\eqref{eq:expeq} reads in a forward way on $\R_+$
\[
\Delta x_{i+1} = \theta_2(\Delta x_i),\; i=2,\ldots,N-1 \quad \mbox{ with} \quad \theta_2 = \Phi_2^{-1}\circ \check{\Phi}_2.
\]
  
As $\Phi_2'(x) = x-\Phi_2(x)$, one checks that 
$$
\theta_2'=\frac{\check\Phi_2'}{\Phi_2'(\Phi_2^{-1}(\check\Phi_2))}=\frac{Id_{|\R_+} +\check\Phi_2}{\theta_2-\check\Phi_2}
$$ 
so that $\theta_2$ satisfies the ordinary differential equation ({\em ODE})
\[
\tfrac 12 (\theta_2^2)' = Id_{|\R_+} + \check \Phi_2\cdot(1+\theta_2').
\]

At this stage, noting that, for every $x\!\in \R_+$,  
\[
\check \Phi_2(x) = e^x-1-x = \sum_{k\ge 2}\frac{x^k}{k!},
\]
we aim at solving this {\em ODE}  by power series i.e. we assume that
\[
\theta_2(x) = \sum_{k\ge 1} a_k x^k
\]
since $\theta_2(0)=\Phi_2^{-1}\circ \check{\Phi}_2(0)= 0$. By standard arguments, we see that 
\[
\tfrac 12 (\theta_2^2)' (x) = \sum_{k\ge 1}b_k x^k \quad\mbox{ with }\quad b_k = \tfrac{k+1}{2} \sum_{\ell=1}^k a_{k+1-\ell}a_{\ell}
\]
and
\[
x+ \check \Phi_2(x)\big(1+\theta_2'(x)\big)= \sum_{k\ge 1} c_k x^k  \quad\mbox{ with }\quad c_k = \Big(\sum_{\ell=0}^{k-2}\frac{\ell+1}{(k-\ell)!}a_{\ell+1} +\frac{1}{k!}\Big) \mbox{\bf 1}_{\{k\ge 2\}} + \mbox{\bf 1}_{\{k=1\}}.
\]

One derives that $a^2_1 = 1$ so that $a_1=1$ since $\theta_2$ is non-decreasing, $\frac{a_1+1}{2}=3 a_1a_2$ which implies $a_2= \frac 13$ and 
\[
a_k = \frac{1}{k+1}\left(\frac{2}{k!} + \sum_{\ell=2}^{k-1} \frac{\ell}{(k+1-\ell)!}a_{\ell}\right) -\frac 12 \sum_{\ell=2}^{k-1} a_{\ell}a_{k+1-\ell}, \quad k\ge 3.
\] 

As a consequence, we can compute $\theta_2(x)$ with an arbitrary accuracy.  Then, the master equation reduces to  the scalar boundary condition
\[
\sum_{k=0}^{N-2}\theta_2^{\circ k}(\Delta x_2)= b
\]
which can be solved numerically by various elementary zero search methods like dichotomy, Newton-Raphson algorithm, etc, since for $k\ge 1$ the $k$-fold composition $\theta_2^{\circ k}$ of $\theta_2$ is continuous and increasing on $\R_+$ from $\theta_2^{\circ k}(0)=0$ to $\theta_2^{\circ k}(+\infty)=+\infty$. Then $x^{N,1,0,b}_k=\sum_{j=0}^{k-2}\theta_2^{\circ j}(\Delta x_2)$ for $k=1,\ldots,N$.

\medskip
\noindent {\bf Numerical example.} The optimal quadratic dual $11$-quantizer of the truncated exponential distribution with parameter $\lambda = 1$ over the unit interval ($a=0$, $b=1$) is given by
\[
x^{11,1,0,1}=(0,0.086271,       0.17510,      0.26663,      0.36105,      0.45853,      0.55929,      0.66355,      0.77156,      0.88361,1).
\]
CPU time on  a $1.8$~MHz processor with Matlab: $6.10^{-3}\,s$ (using a dichotomy algorithm to determine $\Delta x_2$). 

\medskip
\noindent -- {\em General case} ($r\!\in (1, +\infty)$). Set $\check \Phi_r(x) = \Phi_r(-x)$  for every $x\ge 0$. 
By an obvious change of variable, one has 
\[
\check \Phi_r(x)= x^r \int_0^1 \Psi_r(1-z) e^{zx}dz, \quad x\ge 0,
\]
so that $\check \Phi_r(0)=0$, $\check \Phi_r$ is increasing on $\R_+$, goes to infinity at infinity by the monotone convergence theorem since $\Psi_r >0$ on $(0,1)$. As a consequence it is   a $C^1$ homeomorphism of $(0, +\infty)$ (in fact a diffeomorphism since $\check \Phi'_r$ is never $0$ on $(0, +\infty)$). Equation~\eqref{eq:expeq} can be written in a backward way
\[
\Delta x_i = \widetilde \theta_r (\Delta x_{i+1}), \; i= 2,\ldots, N-1\quad \mbox{ with } \quad \widetilde \theta_r = (\check \Phi_r)^{-1}\circ \Phi_r.
\]
Now let us focus for a while  on $\Phi_r$ itself on $\R_+$. As $\Psi_r(z)\ge z^{r+1} $ on $[0,1]$ (see~\eqref{eq:Psi_r}), one has for every $x>0$, 
\begin{align*}
\Phi_r(x) &\ge x^r e^{-x} \int_0^1 z^{r+1}e^{x z}dz = x^{-2} e^{-x} \int_0^x z^{r+1}e^z dz \\
& = x^{r-1} -(r+1)x^{-2}e^{-x}\int_0^x z^re^z dz\\
&\ge x^{r-1}-(r+1) x^{r-2}(1-e^{-x}) \ge \tfrac 12 x^{r-1} \quad \mbox{ for }\quad x\ge 2(r+1).
\end{align*}
Hence $\lim_{x\to +\infty} \Phi_r(x) = +\infty$ which in turn implies that $\lim_{x\to +\infty} \widetilde \theta_r (x)  = +\infty$. Consequently the continuous function $x\mapsto \sum_{k=1}^{N-1}\widetilde \theta_r^{\circ (N-k)}(x)$ is null at $0$ and goes to infinity at infinity so that the equation
\[
\sum_{k=0}^{N-2}\widetilde  \theta_r^{\,\circ k}(x) = b
\]
always has a solution $x^{r,b}$ and we may set $\Delta x_{i} = \widetilde \theta^{\circ (N-i)}_r(x^{r,b})$, $i=2,\ldots,N$. We know that the solution is unique by Theorem~\ref{thm:unique}. 

Unfortunately, we have no semi-closed form for $\widetilde \theta_r$ like in the quadratic case for $\theta_2$  since we could  not find an {\em ODE} satisfied by $\widetilde \theta_r$ in full generality.  

When $r$ is an integer, then $\Psi_r$ is a polynomial function with degree $r-1$ if $r$ is even and $r+1$ if $r$ is odd, whose coefficients of degrees $0$ and $r$ are  always $0$. Then, having in mind that 
\[
\forall\, n\!\in \N, \quad e^{-x} \int_0^x z^ne^{z}dz = (-1)^n n! \Big( \sum_{k=0}^n (-1)^k \frac{x^k}{k!}-e^{-x} \Big)
\] 
it follows that $\Phi_r$ reads
\[
\Phi_r(x) = {\rm sign}(x)^r\left(P_r(x)-e^{-x} \Big(\frac{c_r}{x^2}+Q_r(x)\Big) \right)
\]
where $P_r$ and $Q_r$ are  polynomial functions with degree $r-1$ and $r-2$ respectively that can be computed explicitly and   $c_r=0$ if $r$ is even. 

\medskip 
\noindent -- Case $r=1$. When $r=1$, $\theta_1(z)= 2z(1-z)$ so that $\Psi_1(z)  = 2z^2$. Hence
\begin{align}
\nonumber\Phi_1(x) &= \frac{2\,{\rm sign}(x)}{x^2} e^{-x}\int_0^x z^2e^{z}dz = \frac {4 \,{\rm sign}(x)}{x^2}\Big(\frac{x^2}{2}-x+1-e^{-x}\Big),\\
\label{eq:checkPhi1} \check\Phi_1(x) & =  \frac {4 \,{\rm sign}(x)}{x^2}\Big(e^x -1-x-\frac{x^2}{2}\Big).
\end{align}
In particular,  $\Phi_1$ is increasing on $\R_+$, $C^1$, with $ \lim_{x\to +\infty} \Phi_1(x) = 2$ so that $ \Phi_1$ is  a $C^1$-diffeomorphism  from $(0, +\infty)$ to $(0,2)$. Moreover, it is clear that  $\check \Phi_1$ is a $C^1$-diffeomorphism of $(0, +\infty)$.

In that case we can again write the equation in a forward way 
\[
\Delta x_{i+1} = \theta_1(\Delta x_i),\; i=2,\ldots,N-1 \quad \mbox{ where } \quad \theta_1 = \Phi_1^{-1}\circ \check{\Phi}_1
\]
is defined, $C^1$, increasing  on $\big[0,(\check \Phi_1)^{-1}(2)\big)$ non-negative, satisfies $\theta_1(0) =0$ and $\lim_{x\to \check \Phi_1^{-1}(2)}\theta_1(x) = +\infty$. Consequently the mapping $x\mapsto \sum_{k=1}^{N-1}  \theta_1^{\circ k}(x)$ is defined on an open bounded interval with left endpoint $0$, null at $0$ and goes to infinity at its right endpoint  so that the equation
\[
\sum_{k=0}^{N-2}  \theta_1^{\,\circ k}(x) = b
\]
always has a solution $x^{1,b}$ and we may set $\Delta x_{i} = \theta_1^{\,\circ (i-2)}(x^{1,b}),\;i=2,\ldots,N$.

Moreover,  $\Phi_1$ satisfies the following {\em ODE} on $\R\setminus \{0\}$,
\[
\Phi'_1(x) = -\Big(\frac 2x+1\Big) \Phi_1(x) + 2\,{\rm sign}(x)
\]
from which we derive that, for $x>0$ small enough,
\[
\theta'_1(x) =\frac{(\check \Phi_1)'}{\Phi_1'\circ \theta_1} (x)=  \frac{\check{\Phi}'_1(x)}{2-\big(\frac{2}{\theta_1(x) }+1  \big) \Phi_1(\theta_1(x)) } =  \frac{\check{\Phi}'_1(x)}{2-\big(\frac{2}{\theta_1(x) }+1  \big) \check \Phi_1(x) } 
\]
which can be rewritten (in a neighbourhood of $0$ on $\R_+$) as  the non-linear  {\em ODE}
\[
2\check \Phi_1(x) \theta_1'(x) - (\theta_1^2)'(x) \big(1- \tfrac 12 \check\Phi_1(x) \big) + \theta_1(x) \check{\Phi}'_1(x) = 0,\quad \theta_1(0)=0.
\]

This {\em ODE} can be solved as a power series with positive convergence radius. First note that 
$$
\forall\, x\ge 0, \quad \check \Phi_1(x) =  \sum_{k\ge 1} b_k x^k \quad \mbox{ and  } \quad b_k = \frac{4}{(k+2)!}
$$ 
(so that $b_1= \tfrac 23$, $b_2= \tfrac 16$, etc) owing to~\eqref{eq:checkPhi1}. Assume {\em a priori} that $\theta_1$ can be expanded as 
\[
\theta_1(x)  = \sum_{k\ge 1} a_k x^k, \quad x\ge 0.
\]
Then, if we set
\[
\widetilde a^{(2)}_k = (k+1) \sum_{\ell=1}^{k} a_{\ell} a_{k+1-\ell},\; k\ge 1,\mbox{ so that }\;(\theta_1^2)'(x)=\sum_{k\ge 1}a^{(2)}_kx^k,
\]
elementary though tedious computations show that  the sequence $(a_k)_{k\ge 1}$ satisfies the following induction formula (with  the convention $\sum_{\varnothing}=0$)
\[
a_1=1,\quad a_k = \frac{3}{2(k+2)}\left(\sum_{\ell=1}^{k-1}\left(\widetilde a^{(2)}_{\ell}\frac{b_{k-\ell}}{2} + (k-\ell+1) a_{\ell}b_{k-\ell+1}\right)-(k+1) \sum_{\ell=2}^{k-1} a_{\ell}a_{k+1-\ell}\right), \; k\ge 2,
\]


\medskip
\noindent {\bf Remark.} Another (less tractable) inductive formula can be derived by dealing  directly  with the identity ${\Phi}_1\circ \theta_1 = \check\Phi_1$.



\small
\begin{thebibliography}{0}

\bibitem{AlJo} {\sc  Alfonsi, A. Corbetta J.  and  Jourdain, B.} (2020). Sampling of probability measures in the convex order by Wasserstein projection, {\em Annales de l'Institut Henri Poincar\'e B, Probabilit\'es et Statistiques}, {\bf 56}(3):1706-1729.

\bibitem{AlJo2} {\sc  Alfonsi, A. Corbetta  J. and  Jourdain, B.} (2019). Sampling of one-dimensional probability measures in the convex order and computation of robust option price bounds, {\em International Journal of Theoretical and Applied Finance}, {\bf 22}(3).
  
%
%
%
 %

\bibitem{DeMarch} {\sc De March, H.} (2018).  Entropic approximation for multi-dimensional martingale optimal transport, arXiv 1812.11104.

%
%
%
%
%
\bibitem{IEEE}  {\sc Gersho A. and Gray R.M. eds} (1982). Special issue on quantization. {\em IEEE Trans. Inform. Theory}, {\bf  28} (2), Vol I \& II.

\bibitem{GrLu} {\sc Graf, S. and  Luschgy, H.} (2000). {\em Foundations of quantization for probability distributions},   LNM~1730, Springer, Berlin, 230p.

\bibitem{GuOb}   {\sc Guo, G. and  Obl\`oj, J.}  (2017). Computational Methods for Martingale Optimal Transport problems, {\em Ann. Appl. Probab.}, {\bf 29 }(6):3311--3347.

\bibitem{HL19}   {\sc Henry-Labord\`ere, P.} (2019). (Martingale) optimal transport and anomaly detection with neural networks : a primal-dual algorithm,  arXiv:1904.04546.
  
  %
%
%
%
%
%
%
\bibitem{JoPa1} {\sc   Jourdain, B. and Pag\`es, G.} (2020).  Quantization and martingale couplings. In progress

\bibitem{JoPa2}  {\sc   Jourdain, B. and Pag\`es, G.} (2019). Convex order, quantization and monotone approximations of ARCH models. {\em 	ArXiv:1910.00799}.

\bibitem{Kieff} {\sc Kieffer, J. C.}  (1982). Exponential rate of convergence for Lloyd's method~I. {\em IEEE Trans. Inform. Theory}, {\bf  28} (2):205-210. 


\bibitem{Lloyd} {\sc Lloyd, S.P.}   (1957).  Least squares quantization in PCM, {\em IEEE Transactions on Information Theory}, {\bf 28}(2):129--137 (reprinted from a Bell Telephone Memorandum labs, Murray Hill, NJ, 1957).

\bibitem{McQueen} {\sc MacQueen, J}. (1967). Some methods for classification and analysis of multivariate observations. In {\em Proceedings of the fifth Berkeley symposium on mathematical statistics and probability}, pp. 281-297.

\bibitem{LaPa} {\sc Lamberton, D. and Pag\`es, G.} (1996). On the critical points of the 1-dimensional Competitive Learning Vector Quantization Algorithm. Proceedings of the ESANNı96, (ed., M. Verleysen), Editions D Facto, Bruxelles, 97--106.


\bibitem{PagSpring2018} {\sc Pag\`es, G.} (2018). {\em Numerical Probability: an introduction with applications to Finance}, Springer-Verlag, xvi +579p.

%

\bibitem{Pag2015} {\sc Pag\`es, G.} (2015). Introduction to optimal quantization for numerics,  {\em ESAIM Proc. \& Surveys}, {\bf 48}:29--79.

\bibitem{PagCvx}  {\sc Pag\`es, G.} (2016). Convex order for path-dependent derivatives: a dynamic programming approach. {\em  S\'eminaire de Probabilit\'es XLVIII},  C. Donati, A. Lejay, A. Rouault eds, LNM 2168, Springer, Cham, 33--96. 

\bibitem{PaWi0}  {\sc Pag\`es, G. and  Wilbertz, B.} (2012). Dual Quantization for random walks with application to credit derivatives, {\em J. Comp. Finance}, {\bf 16}(2):33--60.

\bibitem{PaWi1}  {\sc Pag\`es, G. and  Wilbertz, B.} (2012). Intrinsic stationarity for vector quantization: foundation of dual quantization. {\em SIAM J. Numer. Anal.} {\bf  50}(2):747--780.

\bibitem{PaWi2}  {\sc Pag\`es, G. and   Wilbertz, B.} (2012). Optimal Delaunay and Voronoi quantization schemes for pricing American style options. {\em Numerical methods in Finance}, 171--213, Springer Proc. Math., {\bf 12}, Springer, Heidelberg.

\bibitem{PaWi3} {\sc Pag\`es, G. and  Wilbertz, B.} (2018). Sharp rate for the dual quantization problem, {\em S\'eminaire de Probabilit\'es XLV}, C. Donati, A. Lejay, A. Rouault eds, LNM 2215, Springer,  Cham, 119--164.



\bibitem{Trushkin} {\sc Trushkin A.V.}(1982),   Sufficient conditions for uniqueness of a locally optimal quantizer for a class of convex error weighting functions. {\em IEEE Trans. Inform. Theory}, {\bf 28 }(2):187--198.

\bibitem{Str} {\sc Struwe M.} (1990).{\em Variational Methods (Application to non linear
p.d.e \& Hamiltonian Systems)}, Springer, 244p.  

\end{thebibliography}
\end{document}